\documentclass[12pt,twoside]{article}
\usepackage[all]{xy}
\usepackage {amssymb,latexsym,amsthm,amscd,amsmath}

\newcommand{\Z}{\ensuremath{\mathbb{Z}}}
\newcommand{\Q}{\ensuremath{\mathbb{Q}}}
\newcommand{\tr}{\ensuremath{{}^t\!\!}}
\newcommand{\tra}{\ensuremath{{}^t\!}}
\newcommand{\C}{\mathfrak C}
\newcommand{\D}{\mathfrak D}
\newcommand{\E}{\mathcal E}
\newcommand{\F}{\mathcal F}
\newcommand{\Aut}{\textup{Aut}}
\newcommand{\Hom}{\textup{Hom}} 
\newcommand{\End}{\textup{End}}
\newcommand{\cross}{\times}
\newcommand{\tensor}{\otimes}
\newtheorem{theorem}{Theorem}[section]
\newtheorem{lemma}{Lemma}[section]
\newtheorem{corollary}{Corollary}[section]
\newtheorem{proposition}{Proposition}[section]
\newtheorem*{definition}{Definition}
\newtheorem*{theorem*}{Theorem}
\theoremstyle{definition}

\numberwithin{equation}{subsection}

\newcommand{\ignore}[1]{}

\newcommand{\mynote}[1]{}
\begin{document}
\setcounter{section}{0}
\title{\bf REALITY PROPERTIES OF CONJUGACY CLASSES IN $G_2$ }
\author{BY \\ Anupam Singh AND Maneesh Thakur \\ \small Harish-Chandra Research Institute, Chhatnag Road \\ \small Jhunsi, Allahabad 211019, India \\ \small e-mail: anupam/maneesh@mri.ernet.in}
\date{}
\maketitle   

{\bf \it appeared in Isreal Journal of Mathematics 145(2005), 157-192}
\vskip3mm

\begin{abstract}
Let $G$ be an algebraic group over a field $k$.  
We call $g\in G(k)$ {\bf real} if $g$ is conjugate to $g^{-1}$ in $G(k)$. 
In this paper we study reality for groups of type $G_2$ over fields of characteristic different from $2$. 
Let $G$ be such a group over $k$. 
We discuss reality for both semisimple and unipotent elements.  
We show that a semisimple element in $G(k)$ is real if and only if it is  a product of two involutions in $G(k)$. 
Every unipotent element in $G(k)$ is a product of two involutions in $G(k)$. 
We discuss reality for $G_2$ over special fields and construct examples to show that reality fails for semisimple elements in 
$G_2$ over $\Q$ and $\Q_p$. 
We show that semisimple elements are real for $G_2$ over $k$ with $cd(k)\leq 1$.  
We conclude with examples of nonreal elements in $G_2$ over $k$ finite, with characteristic $k$ not $2$ or $3$, which are not 
semisimple or unipotent. 
\end{abstract}


\section{Introduction}
Let $G$ be an algebraic group over a field $k$.  
It is desirable, from the representation theoretic point of view, to study conjugacy classes of elements in $G$. 
Borrowing the terminology from (\cite{fz}), we call an element $g\in G(k)$ {\bf real} if $g$ is conjugate to $g^{-1}$ in $G(k)$.
An {\bf involution} in $G(k)$ is an element $g\in G(k)$ with $g^2=1$.  
Reality for classical groups over fields of characteristic $\neq 2$ 
has been studied in \cite{mvw} by Moeglin, Vign\'{e}ras and Waldspurger. 
That every element of a symplectic group over fields of characteristic $2$ is a product of two involutions is settled in \cite{ni}. 
Feit and Zuckermann discuss reality for spin groups and symplectic groups in \cite{fz}. 
It is well known that every element of an orthogonal group is a product of two involutions (see~\cite{wa} and~\cite{w2}). 
We plan to pursue this for exceptional groups. 
In this paper, we study this property for groups of type $G_2$ over fields of characteristic different from $2$, for both semisimple and unipotent elements.
By consulting the character table of $G_2$ over finite fields in \cite{cr}, one sees that reality is not true for arbitrary elements of $G_2$ (see also Theorem~\ref{finitecounter} and Theorem~\ref{finitecounter'}, in this paper). 
Let $G$ be a group of type $G_2$ over a field $k$ of characteristic $\neq 2$.
We prove that every unipotent element in $G(k)$ is a product of two involutions in $G(k)$. 
As it turns out, the case of semisimple elements in $G(k)$ is more delicate.  
We prove that a semisimple element in $G(k)$ is real in $G(k)$ if and only if it is a product of two involutions in $G(k)$ (Theorem~\ref{maintheorem}).  
We call a torus in $G$ {\bf indecomposable} if it can not be written as a direct product of two subtori, {\bf decomposable} otherwise. 
We show that semisimple elements in decomposable tori are always real (Theorem~\ref{red}). 
We construct examples of indecomposable tori in $G$ containing non-real elements (Proposition~\ref{exsuh} and Theorem~\ref{exsln}). 
We work with an explicit realization of a group of type $G_2$ as the automorphism group of an octonion algebra. 
It is known (Chap. III, Prop. 5, Corollary, \cite{se}) that for a group $G$ of type $G_2$ over $k$, there exists an octonion algebra $\C$ over $k$, unique up to a $k$-isomorphism, such that $G\cong \Aut(\C)$, the group of $k$-algebra automorphisms of $\C$. 
The group $G$ is $k$-split if and only if the octonion algebra $\C$ is split, otherwise $G$ is anisotropic and $\C$ is necessarily a division algebra. 
We prove that any semisimple element in $G(k)$, either leaves invariant a quaternion subalgebra or fixes a quadratic \'{e}tale subalgebra pointwise (Lemma~\ref{fixsubalgebra}). In the first case, reality is a consequence of a theorem of Wonenburger (Th. 4, [W1]). 
In the latter case, the semisimple element belongs to a subgroup $SU(V,h)\subset G$, 
for a hermitian space $(V,h)$ of rank $3$ over a quadratic field extension $L$ of $k$, or to a subgroup $SL(3)\subset G$.  
We investigate these cases separately in sections 6.1 and 6.2 respectively.
We discuss reality for $G_2$ over special fields (Proposition~\ref{exsuh}, Theorem~\ref{exsln} and Theorem~\ref{finitecounter}). 
We show that for $k$ with $cd(k)\leq 1$ (e.g., $k$ a finite field), every semisimple element in $G(k)$ is a product of two involutions in $G(k)$, and hence is real (Theorem~\ref{cdk}). 
We show that nonreal elements exists in $G_2$ over $k$ finite, 
with characteristic $k$ not $2$ or $3$ (compare with [CR]), these are not semisimple or unipotent.
We include a discussion of conjugacy classes of involutions in $G(k)$ over special fields. 
The work of \cite{mvw} has played an important role in representation theory of $p$-adic groups. 
We hope the results in this paper will find applications in the subject.
\section{The Group $G_2$ and Octonions}\label{octonions}
We begin by a brief introduction to the group $G_2$. 
Most of this material is from \cite{sv}. 
Any group $G$ of type $G_2$ over a given field $k$ can be realized as the group of $k$-automorphisms of an octonion algebra over $k$, determined uniquely by $G$. 
We will need the notion of a composition algebra over a field $k$.
\begin{definition}
A composition algebra $\C$ over a field $k$ is an algebra over $k$, not necessarily associative, with an identity element $1$ together with a nondegenerate quadratic form $N$ on $\C$, permitting composition, i.e., $N(xy) = N(x)N(y)\ \ \forall \ x,y \in \C$.
\end{definition}
The quadratic form $N$ is called the {\bf norm} on $\C$. 
The associated bilinear form $N$ is given by : $N(x,y)=N(x+y)-N(x)-N(y)$. 
Every element $x$ of $\C$ satisfies the equation $x^2-N(x,1)x+N(x)1 = 0.$ 
There is an involution (anti automorphism of order $2$) on  $\C$ defined by $\bar x = N(x,1)1-x$. 
We call $N(x,1)1=x+\overline{x}$, as the {\bf trace} of $x$. 
The possible dimensions of a composition algebra over $k$ are $1, 2, 4, 8$. 
Composition algebras of dimension $1$ or $2$ are commutative and associative, those of dimension $4$ are associative but not commutative (called {\bf quaternion} algebras), and those of dimension $8$ are neither commutative nor associative (called {\bf octonion} algebras).

Let $\C$ be an octonion algebra and $G = \Aut (\C)$ be the automorphism group. 
Since any automorphism of an octonion algebra leaves the norm invariant, $\Aut(\C)$ is a subgroup of the orthogonal group $O(\C,N)$. 
In fact, the automorphism group $G$ is a subgroup of the rotation group $SO(N)$ and is contained in $SO(N_1)=\left\{t\in SO(N) \mid t(1) = 1\right\}$, 
where  $N_1 = N|_{1^\perp}$. 
We have (Th. 2.3.5, \cite{sv}),
\begin{proposition}
The algebraic group $\mathcal G = \Aut(\C_K)$, where $\C_K = \C\tensor K$ and $K$ is an algebraic closure of $k$, is the split, connected, simple algebraic group of type $G_2$. 
Moreover, the automorphism group $\mathcal G$ is defined over $k$. 
\end{proposition}
In fact (Chap. III, Prop. 5, Corollary,~\cite{se}), any simple group of type $G_2$ over a field $k$ is isomorphic to the automorphism group of an octonion algebra $\C$ over $k$. 
There is a dichotomy with respect to the norm of octonion algebras (in general, for composition algebras). 
The norm $N$ is a {\bf Pfister} form (tensor product of norm forms of quadratic extensions) and hence is either anisotropic or 
hyperbolic. 
If $N$ is anisotropic, every nonzero element of $\C$ has an inverse in $\C$. 
We then call $\C$ a {\bf division} octonion algebra. 
If $N$ is hyperbolic, up to isomorphism, there is only one octonion algebra 
with $N$ as its norm,  called the {\bf split} octonion algebra. 
We give below a model for the split octonion algebra over a field $k$. 
Let
$$
\C =  \left \{ \left (\begin{array}{cc} \alpha &v \\  w&\beta \\ \end{array}  \right)  | \alpha ,\beta \in k  ; v,w \in k^3 \right\} ,
$$ 
where $k^3$ is the three-dimensional vector space over $k$ with standard basis. 
On $k^3$ we have a nondegenerate bilinear form, given by 
$\langle v,w\rangle=\sum_{i=1}^{3} v_iw_i$, where $v=(v_1,v_2,v_3)$ and $w=(w_1,w_2,w_3)$ in $k^3$ and the wedge product on $k^3$ 
is given by $\langle u\wedge v,w\rangle =\det (u,v,w)$ for $u,v,w \in k^3$. 
Addition on $\C$ is entry-wise and the multiplication on $\C$ is given by,
$$ 
\left (\begin{array}{cc} \alpha &v \\  w&\beta \\ \end{array}  \right)     \left (\begin{array}{cc} \alpha ' &v' \\  w'&\beta ' \\ \end{array}  \right)  =       \left (\begin{array}{cc} \alpha\alpha ' -\langle v,w'\rangle &\alpha v'+\beta ' v + w\wedge w' \\ \beta w'+\alpha ' w + v\wedge v'&\beta\beta '-\langle w,v'\rangle \\ \end{array}  \right).
$$
The quadratic form $N$, the norm on $\C$, is given by 
$$
N\left (\begin{array}{cc} \alpha &v \\  w&\beta \\ \end{array}  \right)  = \alpha\beta + \langle v,w\rangle.
$$  
An octonion algebra over a field $k$ can be defined as an algebra over $k$ which, after changing base to a separable closure $k_s$ of $k$, becomes isomorphic to the split octonion algebra over $k_s$ (see~\cite{t}). 

\subsection{Octonions from rank $3$ Hermitian spaces}\label{octrank3}
We briefly recall here from \cite{t}, a construction of octonion algebras from rank $3$ hermitian spaces over a quadratic \'{e}tale algebra over $k$.  
First we recall (cf.~\cite{kmrt}),
\begin{definition}
Let $\E$ be a finite dimensional $k$-algebra. 
Then $\E$ is called an \'{e}tale algebra if $\E\tensor_k k_s \cong k_s\cross\ldots \cross k_s$, where $k_s$ is a separable closure of $k$.
\end{definition} 
Let $L$ be a quadratic \'{e}tale algebra over $k$ with $x\mapsto \overline{x}$ as its standard involution. 
Let $(V,h)$ be a rank $3$ hermitian space over $L$, i.e., $V$ is an $L$-module of rank $3$ and $h\colon V\cross V\longrightarrow L$ is a nondegenerate hermitian form, linear in the first variable and sesquilinear in the second. 
Assume that the discriminant of $(V,h)$ is trivial, i.e., $\bigwedge^3(V,h)\cong (L, <1>)$, where $<1>$ denotes the hermitian form $(x,y)\mapsto x\overline{y}$ on $L$. 
Fixing a trivialization $\psi \colon \bigwedge^3(V,h)\cong(L,<1>)$, we define a vector product $\cross\colon V\cross V\longrightarrow V$ by the identity, 
$$
h(u,v\cross w)=\psi(u\wedge v\wedge w),
$$ 
for $u,v,w\in V$. 
Let $\C$ be the $8$-dimensional $k$-vector space $\C=C(L;V,h,\psi) = L\oplus V$. 
We define a multiplication on $\C$ by,
$$
(a,v)(b,w)=(ab-h(v,w),~aw+\overline{b}v+v\cross w),~~a,b\in L,~v,w\in V.
$$ 
With this multiplication, $\C$ is an octonion algebra over $k$ with norm $N(a,v)=N_{L/k}(a)+h(v,v)$. 
Note that $L$ embeds in $\C$ as a composition subalgebra. 
The isomorphism class of $\C$, thus obtained, does not depend on $\psi$. 
One can show that all octonion algebras arise this way. 
We need the following (Th. 2.2, \cite{t}),
\begin{proposition}\label{oh}
Let $(V,h)$ and $(V',h')$ be isometric hermitian spaces with trivial discriminant, over a quadratic \'{e}tale algebra $L$. 
Then the octonion algebras $C(L;V,h)$ and $C(L;V',h')$ are isomorphic, under an isomorphism restricting to the identity map on the subalgebra $L$.
\end{proposition}  
We also need the following
\begin{lemma}\label{ohd}
Let $L$ be a quadratic field extension of $k$. 
Let $(V,h)$ be a rank three hermitian space over $L$ with trivial discriminant. 
For any trivialization $\psi$ of the discriminant, the octonion algebra $\C(L;v,h,\psi)$ is a division algebra, if and only if the $k$-quadratic form on $V$, given by $Q(x)=h(x,x)$, is anisotropic. 
\end{lemma}
We note that a similar construction for quaternion algebras can be done, starting from a rank $3$ quadratic space $V$ over $k$, with trivial discriminant. 
Let $B\colon V\cross V\longrightarrow k$ be a nondegenerate bilinear form. 
Assume that the discriminant of $(V,B)$ is trivial, i.e., $\bigwedge^3(V,B)\cong (k, <1>)$, where $<1>$ denotes the bilinear form $(x,y)\mapsto xy$ on $k$.
 Fixing a trivialization $\psi\colon \bigwedge^3(V,B)\cong(k,<1>)$, we define a vector product $\cross\colon V\cross V\longrightarrow V$ by the identity, $B(u,v\cross w)=\psi(u\wedge v\wedge w),$ for $u,v,w\in V$. 
Let $Q$ be the $4$-dimensional $k$-vector space $Q=Q(k;V,B,\psi) = k\oplus V$. 
We define a multiplication on $Q$ by,
$$
(a,v)(b,w)=(ab-B(v,w),~aw+bv+v\cross w),~~a,b\in k,~v,w\in V.
$$ 
With this multiplication, $Q$ is a quaternion algebra over $k$, with norm $N(a,v)=a^2+B(v,v)$. 
The isomorphism class of $Q$ thus obtained, does not depend on $\psi$. 
One can show that all quaternion algebras arise this way. 
\begin{proposition}\label{qb}
Let $(V,B)$ and $(V',B')$ be isometric quadratic spaces with trivial discriminants, over a field $k$. 
Then the quaternion algebras $Q(k;V,B)$ and $Q(k;V',B')$ are isomorphic.
\end{proposition}  

\section{Some subgroups of $G_2$}\label{galoisgroup}
Let $\C$ be an octonion algebra over a field $k$ of characteristic $\neq 2$. 
Let $L$ be a composition subalgebra of $\C$. 
In this section, we describe subgroups of $G=\Aut(\C)$, consisting of automorphisms leaving $L$ pointwise fixed or invariant. 
We define
\[ 
G(\C/L) = \left\{t\in \Aut (\C) | t(x) = x \ \forall \ x \in L\right\}
\] 
and 
\[ 
G(\C,L) = \left\{t\in \Aut (\C) | t(x) \in L \ \forall \ x \in L\right\}.
\]  
Jacobson studied $G(\C/L)$ in his paper (\cite{j}). 
We mention the description of these subgroups here. 
One knows that two dimensional composition algebras over $k$ are precisely the quadratic \'{e}tale algebras over $k$ (cf. Th. 33.17,~\cite{kmrt}). 
Let $L$ be a two dimensional composition subalgebra of $\C$. 
Then $L$ is either a quadratic field extension of $k$ or $L\cong k\cross k$. 
Let us assume first that $L$ is a quadratic field extension of $k$ and $L=k(\gamma)$, where $\gamma^2=c.1 \neq 0$. 
Then $L^{\perp}$ is a left $L$ vector space via the octonion multiplication. 
Also,
\[
h \colon L^{\perp} \times L^{\perp} \longrightarrow L
\]
\[
h(x,y)= N(x,y)+ \gamma^{-1} N(\gamma x,y),
\]
is a nondegenerate hermitian form on $L^{\perp}$ over $L$. 
Any automorphism $t$ of $\C$, fixing $L$ pointwise, induces an $L$-linear map $t|_{L^{\perp}}\colon L^{\perp}\longrightarrow L^{\perp}$. 
Then we have (Th. 3, \cite{j}), 
\begin{proposition}\label{jacobson1}
Let the notations be as fixed above. 
Let $L$ be a quadratic field extension of $k$ as above. 
Then the subgroup $G(\C/L)$ of $G$ is isomorphic to the unimodular unitary group $SU(L^{\perp},h)$ of the three dimensional space $L^\perp$ over $L$ relative to the hermitian form $h$, via the isomorphism,
\begin{eqnarray*}
\psi \colon G(\C/L) &\longrightarrow& SU(L^{\perp},h) \\
      t & \longmapsto & t|_{L^{\perp}}.
\end{eqnarray*}
\end{proposition}
Now, let us assume that $L$ is a split two dimensional \'{e}tale sublagebra of $\C$. 
Then $\C$ is necessarily split and $L$ contains a nontrivial idempotent $e$. 
There exists a basis $B=\{1,u_1,u_2,u_3,e,w_1,w_2,w_3\}$ of $\C$, called the {\bf Peirce basis} with respect to $e$, such that the subspaces $U=\text{span}\{u_1,u_2,u_3\}$ and $W=\text{span}\{w_1,w_2,w_3\}$ satisfy $U=\{x\in\C\mid ex=0,~xe=x\}$ 
and $W=\{x\in\C\mid xe=0,~ex=x\}$. 
We have, for $\eta\in G(\C/L),~x\in U$,
$$
0=\eta(ex)=\eta(e)\eta(x)=e\eta(x),~~\eta(x)e=\eta(x)\eta(e)=\eta(xe)=\eta(x).
$$ 
Hence $\eta(U)=U$. Similarly, $\eta(W)=W$.
Then we have (Th. 4, \cite{j}), 
\begin{proposition}\label{jacobson2}
Let the notations be as fixed above. 
Let $L$ be a split quadratic \'{e}tale subalgebra of $\C$. 
Then $G(\C/L)$ is isomorphic to the unimodular linear group $SL(U)$, via  the isomorphism given by,
\begin{eqnarray*}
\phi \colon G(\C/L) & \longrightarrow & SL(U) \\
       \eta & \longmapsto & \eta|_{U}.
\end{eqnarray*}
Moreover, if we denote the matrix of $\eta|_U$ by $A$ and that of $\eta|_W$ by $A_1$, with respect to the Peirce basis as above, then 
$\tr A_1 = A^{-1}$.
\end{proposition}
In the model of the split octonion algebra as in the previous section, with respect to the diagonal subalgebra $L$, 
the subspaces $U$ and $W$ are respectively the space of strictly upper triangular and strictly lower triangular matrices. 
The above action is then given by, 
$$
\eta\left (\begin{array}{cc} \alpha &v \\  w&\beta \\ \end{array}  \right) = \left (\begin{array}{cc} \alpha & Av \\  
\tr A^{-1}w&\beta  \end{array}\right).
$$
We now compute the subgroup $G(\C,L)$ of automorphisms of the split octonion algebra, leaving invariant a split quadratic \'{e}tale subalgebra. 
We work with the matrix model for split octonions. 
Up to conjugacy by an automorphism, we may assume that the split subalgebra is the diagonal subalgebra.
We consider the map $\rho$ on $\C$ given by 
\begin{eqnarray*}
\rho \colon \C & \longrightarrow & \C \\
\left (\begin{array}{cc} \alpha &v \\  w&\beta \\ \end{array}  \right)  &\mapsto& \left (\begin{array}{cc} \beta &w \\  v&\alpha \\ \end{array}  \right).  
\end{eqnarray*}
Then $\rho$ leaves the two dimensional subalgebra $L=\left\{\left (\begin{array}{cc} \alpha & 0 \\  0 &\beta \\ \end{array}  \right) | \alpha,\beta \in k \right\} $ invariant and it is an automorphism of $\C$, with $\rho^2=1$.
\begin{proposition}\label{gl2} 
Let $\C$ be the split octonion algebra as above and let $L$ be the diagonal split quadratic \'{e}tale subalgebra. 
Then we have, 
$$
G(\C,L) \cong G(\C/L)\rtimes H,
$$ 
where $H$ is the order two group generated by $\rho$.
\end{proposition}
\noindent {\bf Proof : } Let $h\in G(\C,L)$. 
Then $h|_L=1$ or the nontrivial $k$-automorphism of $L$. 
In the first case, $h\in G(\C/L)$ and in the second, $h\rho \in G(\C/L)$. 
Hence $h=g\rho$ for some $g\in G(\C/L)$. 
Moreover, it is clear that $H$ normalizes $G(\C/L)$ in $\Aut(\C)$.  
Since $H \cap G(\C/L) = \{1\}$, we get the required result. \ \ \ \ $\qed$
 
We now give a general construction of the automorphism $\rho$ of an octonion algebra $\C$, not necessarily split, as above. 
We first recall the Cayley-Dickson Doubling for composition algebras :
\begin{proposition}
Let $\C$ be a composition algebra and $\D\subset \C$ a composition subalgebra, $\D\neq \C$. 
Let $a\in \D^{\perp}$ with $N(a)=-\lambda\neq 0$. 
Then $\D_1=\D\oplus \D a$ is a composition subalgebra of $\C$ of dimension $2 \emph{dim}(\D)$. 
The product on $\D_1$ is given by:
$$
(x+ya)(u+va)=(xu+\lambda\overline{v}y)+(vx+y\overline{u})a,~x,y,u,v\in \D,
$$ 
where $x\mapsto \overline{x}$ is the involution on $\D$. 
The norm on $\D_1$ is given by $N(x+ya)=N(x)-\lambda N(y)$. 
\end{proposition}

Let  $\C$ be an octonion algebra and $L\subset\C$, a quadratic composition subalgebra of $\C$. 
Let $a\in L^{\perp}$ with $N(a)\neq 0$. 
Let $\D=L\oplus La$ be the double as described above. 
Then $\D$ is a quaternion subalgebra of $\C$. 
Define $\rho_1\colon \D\rightarrow \D$ by $\rho_1(x+ya)=\sigma(x)+\sigma(y)a$, where $\sigma$ denotes the nontrivial automorphism of $L$. 
Then $\rho_1$ is an automorphism of $\D$, and clearly $\rho_1^2=1$ and $\rho_1|_L=\sigma$. 
We now repeat this construction with respect to $\D$ and $\rho_1$. 
Write $\C=\D\oplus\D b$ for some $b\in\D^{\perp},~~N(b)\neq 0$. 
Define $\rho\colon \C\rightarrow \C$ by,
$$
\rho(x+yb)=\rho_1(x)+\rho_1(y)b.
$$ 
Then $\rho^2=1$ and $\rho |_L=\sigma$ and $\rho$ is an automorphism of $\C$. 
One can prove that this construction yields the one given above for the split octonion algebra and its diagonal subalgebra. 
We have,
\begin{proposition}\label{gl21}
Let $\C$ be an octonion algebra, possibly division, and $L\subset\C$ a quadratic composition subalgebra. 
Then $G(\C, L)\cong G(\C/L)\rtimes H$, where $H$ is the subgroup generated by $\rho$ and $\rho$ is an automorphism of $\C$ with $\rho^2=1$ and $\rho$ restricted to $L$ is the nontrivial $k$-automorphism of $L$.
\end{proposition}

We mention a few more subgroups of $\Aut(\C)$ before closing this section.
 Let $\D\subset\C$ be a quaternion subalgebra. 
Then we have, by Cayley-Dickson doubling, $\C=\D\oplus \D a$ for some $a\in\D^{\perp}$ with $N(a)\neq 0$. 
Let $\phi\in \Aut(\C)$ be such that $\phi(x)=x$ for all $x\in\D$. 
Then for $z=x+ya\in\C$, we have, $\phi(z)=\phi(x)+\phi(y)\phi(a).$
But $a\in\D^{\perp}$ implies $\phi(a)\in \D^{\perp}=\D a$. 
Therefore $\phi(a)=pa$ for some $p\in \D$ and, by taking norms, we see that 
$p\in SL_1(\D)$. 
In fact, we have (Prop. 2.2.1,~\cite{sv}),
\begin{proposition}\label{glh}
The group of automorphisms of $\C$, leaving $\D$ pointwise fixed, is isomorphic to $SL_1(\D)$, the group of norm $1$ elements of $\D$. 
In the above notation, $G(\C/{\D})\cong SL_1(\D)$.
\end{proposition}

We describe yet another subgroup of $\Aut(\C)$. 
Let $\D$ be as above and $\phi\in \Aut(\D)$. 
We can write $\C=\D\oplus \D a$ as above. 
Define $\widetilde{\phi}\in \Aut(\C)$ by $\widetilde{\phi}(x+ya)=\phi(x)+\phi(y)a$. 
Then one checks easily that $\widetilde{\phi}$ is an automorphism of $\C$ that extends $\phi$ on $\D$. 
These automorphisms form a subgroup of $\Aut(\C)$, which (we will abuse notation and) we continue to denote by $\Aut(\D)$. 

\begin{proposition}
With notations as fixed, we have $G(\C,\D)\cong G(\C/{\D})\rtimes \Aut(\D)$.
\end{proposition}
\noindent{\bf Proof : } Clearly $\Aut(\D)\cap G(\C/{\D})=\{1\}$ and $\Aut(\D)$ normalizes $G(\C/{\D})$. 
Now, for $\psi\in G(\C,\D)$, consider the automorphism $\phi=\psi\widetilde{\psi}^{-1}$. 
Then $\phi$ fixes elements of $H$ pointwise and we have $\psi=\phi\widetilde{\psi}\in G(\C/{\D})\rtimes \Aut(\D)$. \ \ \ \ $\qed$

\section{Involutions in $G_2$}
In this section, we discuss the structure of involutions in $G_2$.
Let $G$ be a group of type $G_2$ over $k$ and $\C$ be an octonion algebra over $k$ with $G=\Aut(\C)$. 
We call an element $g\in G(k)$ an {\bf involution} if $g^2=1$. 
Hence nontrivial involutions in $G(k)$ are precisely the automorphisms of $\C$ of order $2$. 
Let $g$ be an involution in $\Aut(\C)$. 
The eigenspace corresponding to the eigenvalue $1$ of $g\in \Aut(\C)$ is the subalgebra $\D$ of $\C$ of fixed points of $g$ and 
is a quaternion subalgebra of $\C$ (\cite{j}). 
The orthogonal complement $\D^{\perp}$ of $\D$ in $\C$ is the eigenspace corresponding to the eigenvalue $-1$. 
Conversely, the linear automorphism of $\C$, leaving a quaternion subalgebra $\D$ of $\C$ pointwise fixed and, acting as multiplication by $-1$ on $\D^{\perp}$, is an involutorial automorphism of $\C$ (see Proposition~\ref{glh}). 
Let $\rho$ be an involution in $G(k)$ and let $\D$ be the quaternion subalgebra of $\C$, fixed pointwise by $\rho$. 
Let $\rho'=g\rho g^{-1}$ be a conjugate of $\rho$ by an element $g\in G(k)$.
 Then, the quaternion subalgebra $\D'=g(\D)$ of $\C$ is fixed pointwise by $\rho'$.
Conversely, suppose the quaternion subalgebra $\D$ of $\C$ is isomorphic to the quaternion subalgebra $\D'$ of $\C$. 
Then, by a Skolem-Noether type theorem for composition algebras (Cor. 1.7.3, \cite{sv}), there exists an automorphism $g$ of $\C$ such that $g(\D)=\D'$. If $\rho$ denotes the involution leaving $\D$ fixed pointwise, $\rho'=g\rho g^{-1}$ fixes $\D'$ pointwise. 
Therefore, we have,
\begin{proposition}
Let $\C$ be an octonion algebra over $k$. 
Then the conjugacy classes of involutions in $G=\Aut(\C)$ are in bijection with the isomorphism classes of quaternion subalgebras of $\C$.
\end{proposition}  
\begin{corollary}
Assume that $_2Br(k)$, the $2$-torsion in the Brauer group of $k$, is trivial, i.e., 
all quaternion algebras over $k$ are split (for example, $cd(k)\leq 1$ fields). 
Then all involutions in $G(k)$ are conjugates.    
\end{corollary}

We need a refinement of a theorem of Jacobson (Th. 2, \cite{j}), due to Wonenburger (Th. 5, \cite{w1}) and Neumann (\cite{n}), 
\begin{proposition}
Let $\C$ be an octonion algebra over a field $k$ of characteristic different from $2$. 
Then every element of $G$ is a product of $3$ involutions. 
\end{proposition}

We will study in the sequel, the structure of semisimple elements in $G(k)$, in terms of involutions. 
We will show that a semisimple element $g\in G(k)$ is real, i.e., conjugate to $g^{-1}$ in $G(k)$, if and only if $g$ is a product of $2$ involutions in $G(k)$.   

\section{Maximal tori in $SU_n$}\label{sun}
We need an explicit description of maximal tori in the special unitary group of a nondegenerate hermitian space for our work, we discuss it in this section (cf.~\cite{r}, Section 3.4).
 Let $k$ be a field of characteristic different from $2$ and $L$ a quadratic field extension of $k$. 
Let $V$ be a vector space of dimension $n$ over $L$. 
We denote by $k_s$ a separable closure of $k$ containing $L$. 
Let $h$ be a nondegenerate hermitian form on $V$, i.e., $h \colon V\times V \longrightarrow L$ is a nondegenerate $k$-bilinear map such that, 
$$
h(\alpha x,y)=\alpha h(x,y),~~h(x,\beta y)=\sigma(\beta)h(x,y),~~h(x,y)=\sigma(h(y,x)),~\forall ~x,y\in V,~\alpha,\beta\in L,
$$
where $\sigma$ is the nontrivial $k$-automorphism of $L$. 
Let $\E$ be an \'{e}tale algebra over $k$. 
It then follows that the bilinear form $T \colon \E\cross \E \longrightarrow k$, 
induced by the trace $:T(x,y) = tr_{\E/k}(xy)$ for $x,y \in \E$, is nondegenerate.

\begin{lemma}
Let $L$ be a quadratic field extension of $k$. 
Let $\E$ be an \'{e}tale algebra over $k$ containing $L$, equipped with an involution $\sigma$, restricting to the non-trivial $k$-automorphism of $L$. 
Let $\F = \E^{\sigma} = \{x\in \E \mid \sigma(x) = x\}$. 
Let $dim_L(\E)=n$. 
For $u\in \F^*$, define
\begin{eqnarray*}
h^{(u)} & \colon & \E \times \E \longrightarrow L \\
h^{(u)} (x,y)&=& tr_{\E/L}(ux\sigma (y)).
\end{eqnarray*}
Then $h^{(u)} $ is a nondegenerate $\sigma$-hermitian form on $\E$, left invariant by $T_{(\E,\sigma)} = \{\alpha\in \E^* \mid \alpha\sigma (\alpha) = 1\}$, under the action by left multiplication.
\end{lemma}
\noindent {\bf Proof : } That $h^{(u)} $ is a hermitian form is clear. 
To check nondegeneracy, let $ h^{(u)} (x,y)=0 ~\forall y\in \E$. 
Then, $tr_{\E/L}(ux\sigma(y))=0~\forall y \in \E$, i.e., $tr_{\E/L}(xy')=0~\forall y' \in \E$. 
Since $\E$ is \'{e}tale, it follows that $x=0$. 
Therefore $h^{(u)} $ is nondegenerate.
Now let $\alpha\in T_{(\E,\sigma)}$. We have, 
$$
h^{(u)} (\alpha x,\alpha y) =  tr_{\E/L}(u\alpha x\sigma (\alpha y)) = tr_{\E/L}(ux\sigma (y)) = h^{(u)} (x,y).
$$ 
Hence the last assertion.
$\ \ \ \ \ \qed$

\paragraph*{\bf Remark : } We note that $\E = \F \tensor_k L$. 
If we put $\F'=\{x\in \E | \sigma(x) = -x\}$ then $\E = \F\oplus \F'$. 
Further, if $L=k(\gamma)$ with $\gamma^2\in k^*$, then $\F' = \F\gamma$.
\paragraph*{\bf Notation :} In what follows, we shall often deal with situations when, for an algebraic group $G$ defined over $k$, 
and for any extension $K$ of $k$, the group $G(K)$ of $K$-rational points in $G$ coincides with $G(k)\otimes_kK$. 
When no confusion is likely to arise, we shall abuse notation and use $G$ to denote both the algebraic group, as well as its group of 
$k$-points. 
We shall identify $T_{(\E,\sigma)}$ with its image in $U(\E,h^{(u)})$, under the embedding via left homotheties. 
\begin{lemma}
With notations as in the previous lemma, $T_{(\E,\sigma)}$ is a maximal $k$-torus in  $U(\E,h^{(u)} )$, the unitary group of the hermitian space $(\E, h^{(u)})$.
\end{lemma}
\noindent
The proof is a tedious, straight forward computation, we omit it here.
\begin{corollary}
Let $T_{(\E,\sigma)}^1 = \{\alpha\in \E^* | \alpha\sigma (\alpha) = 1, \textup{det}(\alpha) = 1\} $. Then $T_{(\E,\sigma)}^1 \subset SU(\E,h^{(u)} )$ is a maximal $k$-torus.
\end{corollary}

\begin{theorem}\label{torus}
Let $k$ be a field and $L$ a quadratic field extension of $k$. 
We denote by $\sigma$ the nontrivial $k$-automorphism of $L$. 
Let $V$ be a $L$-vector space of dimension $n$ with a nondegenerate $\sigma$-hermitian form $h$. 
Let $T \subset U(V,h)$ be a maximal $k$-torus. Then there exists $\E_T$, an \'{e}tale $L$-algebra of dimension $n$ over $L$, with an involution $\sigma_h$ 
restricting to the nontrivial $k$-automorphism of $L$, such that 
$$
T = T_{(\E_T,\sigma_h)}.
$$ 
Moreover, if $\E_T$ is a field, there exists $u\in \F^*$ such that $(V,h)$ is isomorphic to $(\E_T,h^{(u)} )$ as a hermitian space. 
\end{theorem}
\noindent {\bf Proof :}  Let $A = \End_L(V)$. Then $A$ is a central simple $L$-algebra. 
Let $\E_T = \mathcal Z_A(T)$, the centralizer of $T$ in $A$. 
Note that $T\subset \E_T$. 
The hermitian form $h$ defines the adjoint involution $\sigma_h$ on $A$, 
\begin{eqnarray*}
\sigma_h \colon A &\longrightarrow & A \\
h(\sigma_h(f)(x),y)& =& h(x,f(y)) 
\end{eqnarray*} 
for all $x,y\in V$. 
Then $\sigma_h$ is an involution of second kind over $L/k$ on $A$ (cf.~\cite{kmrt}).
We claim that $\sigma_h$ restricts to $\E_T$: Let  $f\in \E_T$, we need to show $\sigma_h(f)\in \E_T$, i.e., $\sigma_h(f)t = t \sigma_h(f)\ \ \forall \ t \in T$. 
This follows from, 
\begin{eqnarray*} 
h(\sigma_h(f)t(x),y) &=& h(t(x),f(y)) = h(x,t^{-1}f(y)) =h(x,ft^{-1}(y)) \\ 
&=& h(\sigma_h(f)(x),t^{-1}y) = h(t\sigma_h(f)(x),y).
\end{eqnarray*}

We have $T\subset U(V,h)\subset \End_L(V)$ and $\sigma_h$ is an involution on $\End_L(V)$, restricting to the nontrivial $k$-automorphism of $L$. 
There is a canonical isomorphism of algebras with involutions (Chap. I, Prop. 2.15,~\cite{kmrt}), 
$$
(\End_L(V)\tensor_k k_s, \sigma_h )\cong(\End_{k_s}(V) \times \End_{k_s}(V), \epsilon), 
$$ 
where $\epsilon(A,B) = (B,A)$. 
Since $U(V,h) = \{A\in \End_L(V) \mid A\sigma_h(A) = 1\}$, we have, 
$$ 
U(V,h)\tensor_k k_s \cong \{(A,B) \in \End_L(V)\tensor_k k_s \mid (A,B).\epsilon(A,B) = 1\}
$$
$$ 
= \{(A,A^{-1}) \mid A \in \End_{k_s}(V)\}.
$$ 
We thus have an embedding 
$$
T\tensor_k k_s\longrightarrow\End_{k_s}(V) \times \End_{k_s}(V),~A\mapsto (A,A^{-1}).
$$ 
To prove $\E_T$ is \'{e}tale, we may conjugate $T\tensor k_s$ to the diagonal torus in $GL_n(k_s)$. 
The embedding then becomes, 
$$
T\tensor_k k_s \cong (k_s^*)^n\longrightarrow M_n(k_s)\cross M_n(k_s),
$$
$$
(t_1 \ldots t_n)\mapsto (diag(t_1,\ldots, t_n), diag(t_1^{-1},\ldots ,t_n^{-1})).
$$
Now, we have,
$$
\E_T\tensor_k k_s = \mathcal Z_A(T)\tensor_k k_s = \mathcal Z_{A\tensor_k k_s}\left(T\tensor_k k_s\right)
$$ 
$$
\cong\mathcal Z_{M_{n}(k_s) \times M_{n}(k_s)} \left(\{(diag(t_1,\ldots, t_n), diag(t_1^{-1},\ldots ,t_n^{-1})) \mid t_i\in k_s^*\}\right) = k_s^{2n}.
$$
Hence $\E_T$ is an \'{e}tale algebra of $k$-dimension $2n$ and $L$-dimension $n$. 
We have,  $T\subset T_{(\E_T,\sigma_h)}$ and, by dimension count, $T=T_{(\E_T,\sigma_h)}$.
We have on $V$, the natural left $\End_L(V)$-module structure. 
Since $\E_T$ is a subalgebra of $\End_L(V)$ and a field, $V$ is a left 
$\E_T$-vector space of dimension $1$. 
Let $V=\E_T.v$ for $v\neq 0$. 
Let us consider the dual $V^*=\Hom_L (V,L)$, which is a left-$\E_T$-vector space of dimension $1$ via the action: $(\alpha.f)(x)=f(\alpha(x)),~\alpha\in\E_T,~x\in V$.  
We consider the following elements in $V^*$:

\begin{eqnarray*}
\phi_1 & \colon & V=\E_T.v \longrightarrow L \\
 && fv \mapsto h(f(v),v) \\
\phi_2 & \colon & V=\E_T.v \longrightarrow L \\
&& fv \mapsto tr(f).
\end{eqnarray*}
Since $\E_T$ is separable, both these are nonzero elements of $V^*$. 
Hence there exists $u\in \E_T^*$ such that $h(f(v),v) = tr(uf) \forall f\in \E_T$.
We have, 
$$
h(f.v,g.v) = h(f(v),g(v)) = h(\sigma_h(g)f(v),v) = tr(u\sigma_h(g)f) \forall f,g \in \E_T.
$$ 
This will prove the lemma provided we show $u\in F$. 
For any $f\in \E_T$ we have,
$$
tr(\sigma_h(u)f)=tr(\sigma_h(u).\sigma_h(\sigma_h(f)))=\sigma_h(tr(u\sigma_h(f)))
$$
$$
=\sigma_h(h(\sigma_h(f)(v),v))=h(v,\sigma_h(f)(v))=h(f(v),v)=tr(uf).
$$
Since $\E_T$ is separable, the trace form is nondegenerate and hence $\sigma_h(u)=u$. 
The map 
$$
\Phi \colon (V,h)\longrightarrow (\E_T,h^{(u)}),~~fv\mapsto f
$$ is an isometry: 
$$
h^{(u)}(\Phi(fv),\Phi(gv))=tr(u\sigma_h(g)f)=h(fv,gv)
$$ by the computation done above. \ \ \ \ \ $\qed$

\begin{corollary}\label{torussuv}
Let the notations be as fixed above. 
Let $T$ be a maximal torus in $SU(V,h)$. 
Then there exists an \'{e}tale algebra $\E_T$ over $L$ of dimension $n$, such that $T \cong T^1_{(\E_T,\sigma_h)}$.
\end{corollary}
\noindent
{\bf Remark :} The hypothesis in the last assertion in Theorem~\ref{torus}, that $\E_T$ be a field, is only a simplifying assumption. 
The result holds good even when $\E_T$ is not a field.

 Let $T\subset SU(V,h)$ be a maximal torus. 
Then from the proof of Theorem~\ref{torus} we get $\E_T = \mathcal Z_{\End(V)}(T')$ is an \'{e}tale algebra with involution $\sigma_h$ such that $T=T^1_{(\E_T,\sigma_h)}$, here $T'$ is a maximal torus in $U(V,h)$.

\begin{lemma}
With notations as above, $V$ is an irreducible representation of $T$ if and only if $\E_T$ is a field.
\end{lemma}
\noindent {\bf Proof :} Suppose $\E_T$ is not a field. 
Then $\exists 0\neq f\in \E_T$ such that $V\neq ker(f)\neq 0$. 
Put $W=ker(f)\subset V$, which is a $L$-vector subspace. 
We claim that $W$ is a $T$ invariant subspace. 
Let $x\in W, t\in T$.
$$
f(x)=0 \Rightarrow t(f(x)) =0 \Rightarrow f(t(x)) =0\Rightarrow t(x)\in W.
$$
Hence, $T(W)=W$.

Conversely, let $\E_T$ be a field and $0\neq W\subset V$ be a $T$-invariant $L$-subspace of $V$. We shall show that $V=W$.   
We know that $V$ is a one dimensional $\E_T$ vector space. 
Thus, it suffices to show that $W$ is an $\E_T$ subspace of $V$. Suppose first that $k$ is infinite.  
Let $t\in T(k)$ be a regular element (see~\cite{bo}, Prop. 8.8 and the Remark on Page 116). 
Then $\E_T=L[t]$ and we have, for $f(t)\in \E_T,~f(t)(W)=W$, since $W$ is $T$-invariant. Now let $k$ be finite. 
Then $\E_T$ is a finite field and its multiplicative group $\E_T^*$ is cyclic. 
The group $T(k)$, being a subgroup of $\E_T^*$, is cyclic. 
Then a cyclic generator $t$ of $T(k)$ is a regular element and arguing as above, we are done in this case too.
\ \ \ \ \ $\qed$

We defined the notion of {\bf indecomposable tori} in the introduction, these are tori which can not be written as a direct product of subtori. 
\begin{corollary}\label{ind}
Let $T$ be a maximal torus in $SU(V,h)$. 
Then $T$ is indecomposable if and only if $V$ is an irreducible representation of $T$. 
That is if and only if $\E_T$ is a field. 
\end{corollary}
\noindent
{\bf Proof :} By the above lemma, if $V$ is reducible as a representation of $T$, $\E_T$ is not a field. 
Hence it must be a product of at least two (separable) field extensions of $L$, say $\E_T=E_1\cross\ldots\cross E_r$. 
Then from Corollary~\ref{torussuv}, $T=T^1_{\E_T}=T^1_{E_1}\cross\ldots\cross T^1_{E_r}$. 
Hence $T$ is decomposable. 
Conversely, suppose $V$ is irreducible as a representation of $T$. 
Then, by the above lemma, $\E_T$ is a field.
Suppose the torus $T$ decomposes as $T=T_1\cross T_2$ into a direct product of two proper subtori. 
Suppose first that $k$ is infinite.  
Let $t\in T(k)$ be a regular element (see~\cite{bo}, Prop. 8.8 and the Remark on Page 116).  
Then the minimal polynomial ($=$ characteristic polynomial) $\chi(X)$ of $t$ factorizes over $k$, as can be seen by base changing to $k_s$ and conjugating $T$ to the diagonal torus in $SL(n)$. 
Therefore $\E_T=L[X]/\chi(X)$ is not a field, a contradiction. 
Hence $T$ is indecomposable. 
When $k$ is finite, the multiplicative group $\E_T^*$ of $\E_T$ is cyclic and hence $T(k)$ is cyclic. 
A cyclic generator $t$ of $T(k)$ is then regular and we repeat the above argument to reach a contradiction. Hence $T$ is indecomposable.\  \  \  \  $\qed$  

\section{Reality in $G_2$}
Let $G$ be a group of type $G_2$ defined over a field $k$ of characteristic $\neq 2$.
 Then, there exists an octonion algebra $\C$ over $k$ such that $G\cong \Aut(\C)$ (Chap. III, Prop. 5, Corollary,~\cite{se}). 
Let $t_0$ be a semisimple element of $G(k)$. 
We will also denote the image of $t_0$ in $\Aut(\C)$ by $t_0$.  
We write $\C_0$ for the subspace of trace $0$ elements of $\C$.  
In this section, we explore the question if $t_0$ is conjugate to $t_0^{-1}$ in $G(k)$.  
We put $V_{t_0}=ker(t_0-1)^8$. 
Then $V_{t_0}$ is a composition subalgebra of $\C$ with norm as the restriction of the norm on $\C$ (\cite{w1}). 
Let $r_{t_0}=dim(V_{t_0}\cap \C_0)$. 
Then $r_{t_0}$ is $1,3$ or $7$. 
We have, 
\begin{lemma}\label{fixsubalgebra}
Let the notations be as fixed above and let $t_0\in G(k)$ be semisimple. 
Then, either $t_0$ leaves a quaternion subalgebra invariant or fixes a quadratic 
\'{e}tale subalgebra $L$ of $\C$ pointwise. 
In the latter case, $t_0\in SU(V,h)\subset G(k)$ for a rank $3$ hermitian space $V$ over a quadratic field extension $L$ of $k$ or $t_0\in SL(3)\subset G(k)$.
\end{lemma}
\noindent{\bf Proof : } From the above discussion, we see that $r_{t_0}$ is $1,3$ or $7$.
If $r_{t_0}=3$ , $t_0$ leaves a quaternion subalgebra $\D$ of $\C$ invariant. 
As in Proposition \ref{glh}, writing $\C=\D\oplus \D a$ for $a\in \D^{\perp},~N(a)\neq 0$, $t_0$ is explicitly given by $t_0(x+ya)=cxc^{-1}+(pcyc^{-1})a$ for some $c\in\D,~N(c)\neq 0$ and $p\in\D,~N(p)=1$. 
We now assume $r_{t_0}=7$. 
In this case, the minimal polynomial of $t_0$ on $\C_0$ is $(X-1)^7$. 
But since  $t_0$ is semisimple, the minimal polynomial of $t_0$ is a product of distinct linear factors over the algebraic closure. 
Therefore $t_0=1$. 
In the case $r_{t_0}=1$,  $L=V_{t_0}$ is a two dimensional composition subalgebra and has the form $V_{t_0}=k.1\oplus (V_{t_0}\cap \C_0)$, an orthogonal direct sum.
 Let $L\cap \C_0 = k.\gamma$ with $N(\gamma)\neq 0$. 
Since $t_0$ leaves $\C_0$ and $V_{t_0}$ invariant, we have, $t_0(\gamma) = \gamma$ and hence $t_0(x)=x\ \forall x\in L$, 
so that $t_0\in G(\C/L)$. The result now follows from Proposition~\ref{jacobson1} and Proposition~\ref{jacobson2}. 
\ \ \ \ $\qed$

If $t_0$ leaves a quaternion subalgebra invariant, it is a product of  two involutions and hence real in $G(k)$. 
This follows from the following theorem (see Th. 4,~\cite{w1}). 
\begin{theorem}\label{wonen}
Let $\C$ be an octonion algebra. 
If $g$ is an automorphism of $\C$ which maps a quaternion subalgebra $\D$ into itself, then $g$ is a product of two involutory automorphisms.
\end{theorem}
\begin{corollary}
If an automorphism $g$ of $\C$ leaves a nondegenerate plane of $\C_0$ invariant, then it is a product of two involutory automorphisms.
\end{corollary}
We discuss the other cases here, i.e., $t_0$ leaves a quadratic \'{e}tale subalgebra $L$ of $\C$ pointwise fixed.
\begin{enumerate}
\item The fixed subalgebra $L$ is a  quadratic field extension of $k$ and 
\item The fixed subalgebra is split, i.e., $L\cong k\cross k$.
\end{enumerate}

By the discussion in section~\ref{galoisgroup}, in the first case, $t_0$ belongs to $G(\C/L)\cong SU(L^{\perp},h)$ (Proposition~\ref{jacobson1}) and in the 
second case $t_0$ belongs to $G(\C/L)\cong SL(3)$ (Proposition~\ref{jacobson2}). 
We denote the image of $t_0$ by $A$ in both of these cases. 
We analyse further the cases when the characteristic polynomial of $A$ is reducible or irreducible.

\begin{theorem} \label{red}
Let $t_0$ be a semisimple element in $G(k)$ and suppose $t_0$ fixes the quadratic \'{e}tale subalgebra $L$ of $\C$ pointwise. 
Let us denote the image of $t_0$ by $A$ in $SU(L^{\perp},h)$ or in $SL(3)$ 
as the case may be. 
Also assume that the characteristic polynomial of $A$ over $L$ in the first case and over $k$ in the second, is reducible. 
Then $t_0$ is a product of two involutions in $G(k)$.
\end{theorem}
\noindent {\bf Proof : } Let us consider the case when $L$ is a field extension. Let $T$ be a maximal torus in $SU(L^{\perp},h)$ containing $t_0$. 
 By Corollary \ref{torussuv}, there exists an \'{e}tale $L$-algebra $\E_T$ with an involution $\sigma$ and $u\in \F^*$ such that $(L^{\perp},h)\cong (\E_T,h^{(u)})$, here $\F$ is the fixed point subalgebra of $\sigma$ in $\E_T$. 
Since the characteristic polynomial of $A$ is reducible, we see that $L^{\perp}$ is a reducible representation of $T$. 
From Corollary \ref{ind} we see that $\E_T$ is not a field. 
We can write $\E_T\cong \F\tensor L$ where $\F$ is a cubic \'{e}tale $k$-algebra but not a field. 
Let $\F= k\times \Delta$, for some quadratic \'{e}tale $k$-algebra~$\Delta$. 
Hence $\E_T\cong L\times (\Delta\tensor L)$ and $\sigma$ is given by $(\alpha,f\tensor\beta)\mapsto (\bar\alpha,f\tensor\bar\beta) $. 
Writing $u=(u_1,u_2)$ where $u_1\in k$, the hermitian form $h^{(u)}$ is given by $h^{(u)}((l,\delta),(l',\delta'))=tr_{L/L}(lu_1l')+tr_{\Delta\tensor L/L}(\delta u_2\delta')=lu_1l'+tr_{\Delta\tensor L/L}(\delta u_2\delta')$. 
Hence $L\times \{0\}$ is a nondegenerate subspace left invariant by the action of $t_0\in T^1_{(\E_T,\sigma_h)}\cong  
T^1_L\times T^1_{\Delta\otimes L}$, which acts by left multiplication. 
Therefore $t_0$ leaves invariant a two dimensional nondegenerate $k$-plane invariant in $\C_0$. 
The result now follows from Corollary to Theorem \ref{wonen}. 
The proof in the case when $L$ is split proceeds on similar lines. \ \ \ \ $\qed$

In general, we have the following,
\begin{theorem}\label{maintheorem}
Let $G$ be a group of type $G_2$ over a field $k$ of characteristic not $2$.
 Then every unipotent element in $G(k)$ is a product of two involutions in $G(k)$. 
Let $g\in G(k)$ be a semisimple element. 
Then, $g$ is real in $G(k)$ if and only if it is a product of two  involutions in $G(k)$. 
\end{theorem} 
\noindent
{\bf Proof :} The assertion about unipotents in $G(k)$ follows from a theorem of Wonenburger (Th. 4,~\cite{w1}), 
which asserts that if the characteristic polynomial of $t\in \Aut(\C)$ is divisible by $(x-1)^3$, $t$ is a product of two involutory automorphisms of $\C$. \ \ \ \ $\qed$

In view of Theorem~\ref{red}, we need to consider only the semisimple elements in $SU(L^{\perp},h)$ or in $SL(3)$ with irreducible characteristic polynomials. 
By Corollary~\ref{ind}, it follows that such elements lie in indecomposable tori. 
The result follows from the following theorem. 
\begin{theorem} 
Let $t_0$ be an element in $G(k)$ and suppose $t_0$ fixes a quadratic \'{e}tale subalgebra $L$ of $\C$ pointwise. 
Let us denote the image of $t_0$ by $A$ in $SU(L^{\perp},h)$ or in $SL(3)$ 
as the case may be. 
Also assume that the characteristic polynomial of $A$ over $L$ in the first case and over $k$ in the second, is irreducible. 
Then $t_0$ is conjugate to $t_0^{-1}$ in $G(k)$ if and only if $t_0$ is a product of two involutions in $G(k)$.
\end{theorem}
\noindent{\bf Proof :} We distinguish the cases of both these subgroups below and complete the proof in next two subsections, 
see Theorem~\ref{mainsuvh} and Theorem~\ref{counterexample}. \ \ \ \ $\qed$ 

\subsection{$SU(V,h)\subset G$}\label{division}
We assume that $L$ is a quadratic field extension of $k$. 
Let $t_0$ be an element in $G(\C/L)$ with characteristic polynomial of the restriction to $V=L^{\perp}$, irreducible over $L$. 
We write $\C = L\bigoplus V$, where $V$ is an $L$-vector space with hermitian form $h$ induced by the norm on $\C$. 
Then we have seen  that $G(\C/L)\cong SU(V,h)$ (Theorem \ref{jacobson1}). 
\begin{lemma}\label{sufix}
Let the notations be as fixed above. 
Let $t_0$ be an element in $G(\C/L)$ with characteristic polynomial irreducible over $L$. 
Suppose that $\exists g\in G(k)$ such that $gt_0g^{-1}=t_0^{-1}$. 
Then $g(L)=L$.
\end{lemma}
\noindent {\bf Proof : } Suppose that $g(L)\not\subset L$. 
Then we claim that $\exists x \in L\cap \C_0$ such that $g(x)\not\in L$. 
For this, let $y\in L$ be such that $g(y)\not\in L$. 
Let $x=y-\frac{1}{2}tr(y)1$. 
Then $tr(x)=0$ and if $g(x)\in L$ then $g(y)\in L$, a contradiction. 
Hence we have $x\in L\cap \C_0$ with $g(x)\not\in L$. 
Also since $t_0(x)=x$, we have, 
$$
t_0(g(x))=gt_0^{-1}(x)=g(x).
$$ 
Let $g(x)=\alpha+y$, for $0\neq y\in L^{\perp}$, then $t_0(g(x))=\alpha+t_0(y)=\alpha+y$, i.e., $t_0(y)=y$. 
Therefore $t_0$ fixes an element in $L^{\perp}$. 
This implies that the characteristic polynomial of $t_0$ on $L^{\perp}=V$  is reducible, a contradiction.
 Hence, $g(L)=L$. 
\  \  \  \ $\qed$

We recall a construction from Proposition \ref{gl21}. 
Let $a\in L^{\perp}$ with $N(a)\neq 0$. 
Let $\D=L\oplus La$ and $\rho_1 \colon \D \rightarrow\D$ be defined by $\rho_1(x+ya)=\sigma(x)+\sigma(y)a$. 
Write again $\C=\D\oplus \D b$, for $b\in \D^{\perp}$ with $N(b)\neq 0$ and define $\rho\colon \C \rightarrow \C$ by $\rho(x+yb)=\sigma(x)+\sigma(y)b$. 
Then $\rho$ is an automorphism of $\C$ of order $2$ which restricts to $L$ to the nontrivial automorphism of $L$. 
The basis 
$$
\{f_1=a,f_2=b,f_3=ab\}
$$ 
of $V=L^{\perp}$ over $L$ is an orthogonal basis for $h$. 
{\bf We fix this basis throughout this section}. 
Let us denote the matrix of $h$  with respect to this basis by $H=diag(\lambda_1,\lambda_2,\lambda_3)$ where $\lambda_i=h(f_i,f_i)\in k^*$. 
Then $SU(V,h)$ is isomorphic to $SU(H)=\{A\in SL(3,L) \mid \tr AH\bar A=H\}$.
\begin{theorem}\label{sucong2}
With notations fixed as above, let $A$ be the matrix of $t_0$ in $SU(H)$ with respect to the fixed basis described above. 
Let the characteristic polynomial of $A$ be irreducible over $L$. 
Then $t_0$ is conjugate to $t_0^{-1}$ in  $G(k)$, if and only if $\bar A$ is conjugate to $A^{-1}$ in $SU(H)$, where the entries of $\bar A$ are obtained by applying $\sigma$ on the entries of $A$. 
\end{theorem}
\noindent {\bf Proof : }
Let $g\in G(k)$ be such that $gt_0g^{-1}=t_0^{-1}$. 
In view of Lemma~\ref{sufix}, we have $g(L)=L$. 
We have (Prop.~\ref{gl21}) $G(\C, L)\cong G(\C/L)\rtimes N$ where $N=<\rho>$ and $\rho$ is an automorphism of $\C$, described above.  
Clearly $g$ does not belong to $G(\C/L)$. For if so, we can conjugate $t_0$ to $t_0^{-1}$ in $G(\C/L)\cong SU(H)$. 
But then the characteristic polynomial $\chi(X)=X^3-\bar aX^2+aX-1$, where $a\in L$, and $\bar a=a$. Hence $\chi(X)$ is reducible, a contradiction. 
We write $g=g'\rho$ where $g'\in G(\C/L)$. 
Let $B$ be the matrix of $g'$ in $SU(H)$. 
Then, by a direct computation, it follows that,
$$ 
gt_0g^{-1}(\alpha_0.1+\alpha_1 f_1+\alpha_2 f_2+\alpha_3 f_3)
$$
$$
= \alpha_0.1+\alpha_1 B\bar{A}B^{-1}f_1+\alpha_2 B\bar{A}B^{-1}f_2+\alpha_3 B\bar{A}B^{-1}f_3.
$$
Also,
$$
t_0^{-1}(\alpha_0.1+\alpha_1 f_1+\alpha_2 f_2+\alpha_3 f_3)=(\alpha_0.1+\alpha_1 A^{-1}f_1+\alpha_2 A^{-1}f_2+\alpha_3 A^{-1}f_3).
$$
Therefore, if $t_0$ is conjugate to $t_0^{-1}$ in  $G=\Aut(\C)$, then $\bar A$ is conjugate to $A^{-1}$ in $SU(H)$. 
Conversely, let $B\bar AB^{-1}=A^{-1}$ for some $B\in SU(H)$. 
Let $g'\in G(\C/L)$ be the element corresponding to $B$. 
Then $g'\rho$ conjugates $t_0$ to $t_0^{-1}$. \ \ \ \ \ $\qed$

Let $V$ be a vector space over $L$ of dimension $n$ with a nondegenerate hermitian form $h$. 
Let $H$ denote the diagonal matrix of $h$ with respect to some fixed orthogonal basis. 
Then, for any $A\in U(H)$, we have $\tra AH\bar A=H$. 
Let $A\in SU(H)$ with characteristic polynomial $\chi_A(X)=X^n+a_1X^{n-1}+\ldots +a_{n-1}X+(-1)^n$. Then $(-1)^na_i=\bar a_{n-i}$ for $i=1,\ldots,n-1$.
\begin{lemma}
With notations as above, let $A\in SU(H)$ with its characteristic polynomial over $L$ be the same as its minimal polynomial. 
Suppose $A=A_1A_2$ with $A_1,A_2\in GL(n,L)$ and $\bar A_1 A_1=I=\bar A_2 A_2$. 
Then, $A_1,A_2 \in U(H)$.
\end{lemma}
\noindent{\bf Proof : } Let $H=diag(\lambda_1,\lambda_2,\ldots,\lambda_n)$, where $\lambda_1,\ldots,\lambda_n \in k$. We have $\tr AH\bar A=H$. 
Then,
$$
(HA_1^{-1})A(HA_1^{-1})^{-1}=HA_1^{-1}A_1A_2A_1H^{-1}=H\bar A^{-1}H^{-1}=\tr A.
$$
Since the characteristic polynomial of $A$ equals its minimal polynomial, by (Th.~2,~\cite{tz})  $HA_1^{-1}$ is symmetric, i.e., $HA_1^{-1}={}^t(HA_1^{-1})=\tr A_1^{-1}H$. 
This implies, $H=\tr A_1 HA_1^{-1}=\tr A_1H\bar A_1$. 
Hence, $A_1\in U(H)$. 
By similar analysis we see that $A_2\in U(H)$.\ \ \ \ \ $\qed$
\begin{lemma}
With notations as above, let $A\in SU(H)$ with characteristic polynomial $\chi_A(X)=X^n+a_1X^{n-1}+\ldots +a_{n-1}X+(-1)^n$ over $L$, equal to its minimal polynomial. 
Then, $A=B_1B_2$ with $B_1,B_2\in GL(n,L)$ and $\overline B_1 B_1=I=\overline B_2 B_2$. 
\end{lemma}
\noindent{\bf Proof : } Let $A_{\chi}$ denote the companion matrix of $A$, namely 
$$ 
A_{\chi} = \left (\begin{array}{ccccc} 0 & 0 &\ldots&0& -(-1)^n \\ 1 & 0&\ldots&0 &-a_{n-1} \\ \vdots &\vdots&&&\vdots\\ 0&0&\ldots&1 &-a_1 \end{array}  \right) . 
$$
We have, 
$$ 
A_{\chi} = \left (\begin{array}{ccccc} (-1)^n & 0&\ldots&0 & 0 \\ a_{n-1} & 0 &\ldots&0&-1 \\ \vdots&\vdots&&&\vdots\\a_1&-1&\ldots&0 &0 \end{array}  \right) 
\left (\begin{array}{ccccc} 0&0&\ldots & 0 & -1 \\ 0&0&\ldots & -1 &0 \\ \vdots&\vdots&&&\vdots\\ -1&0&\ldots&0 &0\end{array}  \right)
=A_1A_2, 
$$
and $\bar A_1A_1=I=\bar A_2A_2$, using $(-1)^na_i=\bar a_{n-i}$ for $i=1,\ldots,n-1$. 
Since the characteristic polynomial of $A$ equals its minimal polynomial, there exists $T\in GL(n,L)$ such that $A=TA_{\chi}T^{-1}$. 
We put $B_1=TA_1\bar T^{-1}, B_2=\bar TA_2T^{-1}$. 
Then $A=B_1B_2$, where $ \bar B_1B_1=I=\bar B_2B_2$. \ \ \ \ \ $\qed$

\begin{corollary}
Let $A\in SU(H)$ with characteristic polynomial $\chi_A(X)$ over $L$ same as its minimal polynomial. 
Then, $A=B_1B_2$ with $B_1,B_2\in U(H)$ and $\overline B_1 B_1=I=\overline B_2 B_2$. 
\end{corollary}
From this corollary, we get the following,
\begin{lemma}\label{suco1} 
Let $A\in SU(H)$, with characteristic polynomial over $L$ equal to its minimal polynomial. 
Then,
\begin{enumerate}
\item $\bar A$ is conjugate to $A^{-1}$ in  $U(H)$, if and only if $A=A_1A_2$ with $A_1,A_2\in U(H) $ and  $\bar A_1 A_1=I=\bar A_2 A_2$.
\item $\bar A$ is conjugate to $A^{-1}$ in  $SU(H)$, if and only if $A=A_1A_2$ with $A_1,A_2\in SU(H) $ and $\bar A_1 A_1=I=\bar A_2 A_2$.
\end{enumerate}
\end{lemma}

The following proposition is due to Neumann (\cite{n} Lemma 5). 
Recall that we have fixed a basis $\{f_1,f_2,f_3\}$ for $V=L^{\perp}$ over $L$ in Theorem~\ref{sucong2}.
\begin{proposition}\label{neum}
Let $\C$ be an octonion algebra over $k$ and let $L$ be a quadratic field extension of $k$, which is a subalgebra of $\C$. 
An element $t\in G(\C/L)$ is a product of two involutions in $\Aut(\C)$, if and only if, the corresponding matrix $A\in SU(H)$ is a product of two matrices $A_1,A_2\in SU(H)$, satisfying $\bar A_1A_1=\bar A_2A_2=I$.
\end{proposition}
We now have,
\begin{theorem}\label{mainsuvh}
Let $t_0$ be an element in $G(\C/L)$ and let $A$ denote the image of $t_0$ in $SU(H)$. 
Suppose the characteristic polynomial of $A$ is irreducible over $L$. 
Then $t_0$ is conjugate to $t_0^{-1}$, if and only if $t_0$ is a product of two involutions in $G(k)$. 
\end{theorem}
\noindent{\bf Proof : } From Theorem \ref{sucong2} we have, $t_0$ is conjugate to $t_0^{-1}$, if and only if $\bar A$ is conjugate to $A^{-1}$ in $SU(H)$. 
From Lemma \ref{suco1} above, $\bar A$ is conjugate to $A^{-1}$ in  $SU(H)$ if and only if $A=A_1A_2$ with $A_1,A_2\in SU(H) $ and $\bar A_1 A_1=I=\bar A_2 A_2$. 
Now, from Proposition \ref{neum}, it follows that $t_0$ is a product of two involutions.
\ \ \ \ \ $\qed$

Let $V$ be a vector space over $L$ of dimension $n$ together with a nondegenerate hermitian form~$h$. 
Let $A \in SU(H)$.
 Let us denote the conjugacy class of $A$ in $U(H)$ by $C$ and the centralizer of $A$ in $U(H)$ by $\mathcal Z$ and let 
$$
L_A=\{\det(X) \mid X\in \mathcal Z\}.
$$

\begin{lemma}
With notations as fixed above, for $X,Y\in U(H)$, $XAX^{-1}$ is conjugate to $YAY^{-1}$ in $SU(H)$ if and only if $\det(X)\equiv \det(Y)(mod L_A)$.
\end{lemma}
\noindent{\bf Proof : } Suppose there exists $S\in SU(H)$ such that $SXAX^{-1}S^{-1}=YAY^{-1}$. 
Then, $Y^{-1}SX\in \mathcal Z$ and $\det(X)\equiv \det(Y)(mod L_A)$.

Conversely, let $\det(XY^{-1})=\det(B)$ for $B\in \mathcal Z$. 
Put $S=YBX^{-1}$. Then $\det(S)=1$, $S\in SU(H)$ and $Y^{-1}SX=B\in \mathcal Z$. 
Then, $Y^{-1}SXA=AY^{-1}SX$ gives $SXAX^{-1}S^{-1}=YAY^{-1}$. \ \ \ \ \ $\qed$

\begin{lemma}
Let $t_0$ be an element in $G(\C/L)$ for $L$ a quadratic field extension of $k$ and $A$ be the corresponding element in $SU(H)$. 
Suppose the characteristic polynomial of $A$ is irreducible over $L$. 
Then, $t_0$ is conjugate to $t_0^{-1}$ in $G(k)$, if and only if for every $X\in U(H)$ such that $X\bar AX^{-1}=A^{-1}$, $\det(X)\in L_{\bar A}$.
\end{lemma}
\noindent{\bf Proof : }  We have, by Theorem \ref{sucong2}, $t_0$ is conjugate to $t_0^{-1}$ in $G(k)$ if and only if $\bar A$ is conjugate to $A^{-1}$ in $SU(H)$. 
Let $X\in U(H)$ be such that $X\bar AX^{-1}=A^{-1}$. 
Then from the above lemma, $\bar A$ is conjugate to $A^{-1}$ in $SU(H)$ if and only if $\det(X)\in L_{\bar A}$. \ \ \ \ $\qed$
\begin{corollary}\label{suobs}
With notations as fixed above, whenever $L^1/L_{\bar A}$ is trivial, $t_0$ is conjugate to $t_0^{-1}$ in $G(k)$, where $L^1=\{\alpha\in L| \alpha\bar\alpha=1\}$.
\end{corollary}
\noindent{\bf Proof : } We have $L^1=\{\alpha\in L| \alpha\bar\alpha=1\}=\{\det(X)|X\in U(H)\}$. 
Now let us fix $X_0\in U(H)$ such that $X_0\bar AX_0^{-1}=A^{-1}$. 
Then, for any $X\in U(H)$ such that $X\bar AX^{-1}=A^{-1}$, we have $X_0^{-1}X\in \mathcal Z_{U(H)}(\bar A)$. 
Hence $\det(X)\in \det(X_0)L_{\bar A}$. 
But since $L^1/L_{\bar A}$ is trivial, we have $\det(X)\in L_{\bar A}$. 
From the above lemma, it now follows that $t_0$ is conjugate to $t_0^{-1}$ in $G(k)$. \ \ \ \ $\qed$

\noindent{\bf Remark : } From the proof above, for any $X\in U(H)$ such that $X\bar AX^{-1}=A^{-1}$, we get $X\in X_0\mathcal Z_{U(H)}(\bar A)$. 
Since  the characteristic polynomial of $A$ is irreducible, that of $\bar A$ is irreducible as well. 
Therefore $\mathcal Z_{U(H)}(\bar A)\subset \mathcal Z_{\End_L(V)}(\bar A)=L[\bar A]\cong L[T]/<\chi_{\bar A}(T)>$. 
In fact, $\mathcal Z_{U(H)}(\bar A)=\{x\in \mathcal Z_{\End_L(V)}(\bar A) \mid x\sigma_h(x)=1\}$. 
Hence we can write $X=X_0f(\bar A)$ for some polynomial $f(T)\in L[T]$.
\begin{lemma}Let $A\in SU(H)$ and its characteristic polynomial $\chi_A(X)$ be irreducible over $L$. 
Let $\E=L[X]/\chi_{\bar A}(X)$, a degree three field extension of $L$. 
Then $L^1/L_{\bar A}\hookrightarrow L^*/N(\E^*)$.
\end{lemma}
\noindent {\bf Proof : } Define a map $\phi \colon L^1 \longrightarrow L^*/N(\E^*)$ by $x\mapsto xN(\E^*)$. 
We claim that $ker(\phi)=\{x\in L^1\mid x\in N(\E^*)\}=L_{\bar A}=\{N(x)\mid x\in \E^*, x\sigma(x)=1\}$. 
Let $x\in ker(\phi)$, i.e., $x=N(y)$ for some $y\in \E^*$ and $x\sigma(x)=1$. 
Let $\tilde y=xy^{-1}\sigma(y)\in \E^*$ then $N(\tilde y)=x, \tilde y\sigma(\tilde y)=1$. 
Hence $x\in L_{\bar A}$.
 Conversely, if $N(x)\in L_{\bar A}$  for some $x\in \E^*$ such that $x\sigma(x)=1$ then $N(x)\in ker(\phi)$. \ \ \ \ $\qed$

Hence if the field $k$ is $C_1$ (for example, finite field) or it does not admit degree three extensions (real closed fields, algebraically closed fields etc.),  $L^*/N(E^*)$ is trivial. 
From Corollary \ref{suobs}, it follows that every element in $G(\C/L)$, with irreducible characteristic polynomial, is conjugate to its inverse. 
In particular, combining with Theorem \ref{red}, it follows that every semisimple element in $G(k)$ is conjugate to its inverse.

\begin{proposition}\label{counsun}
With notations as above, let $L$ be a quadratic field extension of $k$ and let $S\in SU(H)$ be an element with irreducible characteristic polynomial over $L$, satisfying $\bar S=S^{-1}$. 
Let $\E=L[X]/\chi_S(X)$, a degree three field extension of $L$, and assume $L^1/N(\E^1)$ is nontrivial, where $L^1=\{x\in L\mid x\sigma(x)=1\}$, $\E^1=\{x\in \E\mid x\sigma(x)=1\}$ and $\sigma$ is the extension of the nontrivial automorphism of $L$ to $\E$. 
Then there exists an element $A\in SU(H)$ with characteristic polynomial same as the characteristic polynomial of $S$, which can not be written as $A=A_1A_2$ where $\bar A_i=A_i^{-1}$ and $A_i\in SU(H)$. 
The corresponding element $t$  in $G(\C/L)$ is not a product of two involutions in $G=\Aut(\C)$ and hence not real in $G$.
\end{proposition}
\noindent {\bf Proof : }  Let $b\in L^1$ such that $b^2\not\in N(\E^1)$. 
Put $D=diag(b,1,1)$ and $A=DSD^{-1}$, then $A$ belongs to $SU(H)$. 
Now suppose $A=A_1A_2$ with $\bar A_i=A_i^{-1}$ and $A_i\in SU(H)$. 
Then $A=A_1A_2=DSD^{-1}=DSDD^{-2}$. 
Put $T_1=DSD$ and $T_2=D^{-2}$, then $\bar T_i=T_i^{-1}$. 
Since $A_2AA_2^{-1}=\bar A^{-1}$ and $T_2AT_2^{-1}=\bar A^{-1}$, we have $T_2^{-1}A_2\in \mathcal Z_{U(V,h)}(A)$, i.e., $T_2^{-1}A_2=f(A)$ for some $f(X)\in L[X]$ (see the Remark after Corollary \ref{suobs}). 
Then $b^2=\det(T_2^{-1})=\det(T_2^{-1}A_2)=\det(f(A))\in N(\E^*)$, a contradiction. \ \ \ \ $\qed$

\noindent
{\bf Remark : } If we choose $S$ in the theorem above with characteristic polynomial separable, then the element $A$, constructed in the proof, is a semisimple element in an indecomposable maximal torus, contained in $SU(H)$, which is not real.

We recall that any central division algebra of degree three is cyclic (Theorem, Section 15.6, \cite{p}).  
Let $L$ be a quadratic field extension of $k$. 
Let $F$ be a degree three cyclic extension of $k$ and we denote  $E=F.L$. 
Let us denote the generator of the Galois group of $F$ over $k$ by $\tau$. 
Let $A=F\oplus Fu\oplus Fu^2$ with $udu^{-1}=\tau(d)$ for all $d\in F$ and $u^3=a\in k^*$. 
Then $A$, denoted by $(F,\tau,a)$, is a cyclic algebra of degree three over $k$.  
Recall also that $(F,\tau,a)$ is a division algebra if and only if $a\not\in N_{F/k}(F^*)$. 
We denote the relative Brauer group of $F$ over $k$ by $B(F/k)$, i.e., the group of Brauer classes of central simple algebras over $k$, 
which split over $F$. 
We define a map $\phi\colon B(F/k)\longrightarrow B(E/L)$ by $[(F,\tau,a)]\mapsto [(E,\tau,a)]$ (which is the same as the map $[D]\mapsto [D\otimes L]$). 
This map is well defined (Section 15.1, Cor. c, \cite{p}) and is an injective map since $ker(\phi)=\{[(F,\tau,a)]\in B(F/k)\mid a\in k^*, a\in N_{E/L}(E^*)\}=\{[(F,\tau,a)]\in B(F/k)\mid a\in N_{F/k}(F^*)\}$.  
We have a commutative diagram, 
\[ 
\begin{CD} k^*/N_{F/k}(F^*) @>{\cong}>> B(F/k)\\
@VVV      @VV\phi V \\
L^*/N_{E/L}(E^*) @>{\cong}>> B(E/L)\\
\end{CD} 
\]
The vertical maps are injective in the above diagram. We have the following exact sequence,   
$$
1\longrightarrow (N_{E/L}(E^*)k^*)/N_{E/L}(E^*)\longrightarrow L^*/N_{E/L}(E^*) \longrightarrow L^1/N_{E/L}(E^1) \longrightarrow 1
$$
where $(N_{E/L}(E^*)k^*)/N_{E/L}(E^*) \cong k^*/N_{F/k}(F^*)$. 
Hence, from the commutativity of the above diagram, we get $B(E/L)/\phi(B(F/k)) \cong L^1/N_{E/L}(E^1)$.

This shows $L^1/N_{E/L}(E^1)$ is nontrivial, if and only if there exists a central division algebra $D$ over $L$ which splits over $E$ and it does not come from  a  central division algebra over $k$, split by $F$. 
We recall a proposition from~\cite{k} (Chapter V, Prop. 1). 
\begin{proposition}
Let $k$ be a number field and $L$ a quadratic field extension of $k$. 
Let $F$ be a cyclic extension of degree $m$ over $k$, which is linearly disjoint from $L$, over $k$. 
Then there exists a central division algebra $(FL/L,\tau,a)$ over $L$ of degree $m$, with an involution of second kind, with $a\in L^1$. 
\end{proposition}
\begin{corollary}\label{lne}
Let $k$ be a number field and $L$ a quadratic field extension of $k$. 
Let $F$ be a cyclic extension of degree $3$ over $k$. 
Let us denote $E=F.L$. 
Then $L^1/N_{E/L}(E^1)$ is nontrivial.
\end{corollary}
\noindent{\bf Proof : } By the proposition above, there exists a degree three central division algebra $(E,\tau,a)$ over $L$ with $a\in L^1$. 
Therefore  $a\not\in N_{E/L}(E^*)$. \ \ \ \ $\qed$  

We proceed to construct an example of the situation required in proposition~\ref{counsun}.
\begin{proposition}\label{exsuh}
Let $k$ be a number field. 
There exist octonion algebras $\C$ over $k$ such that not every (semisimple) element in $\Aut(\C)$ is real.
\end{proposition}
\noindent {\bf Proof : } We use Proposition~\ref{counsun} here. 
Let $L$ be a quadratic field extension of $k$. 
Let $F$ be a degree three cyclic extension of $k$. 
Then we have $E=F.L$, a degree three cyclic extension of $L$. 
We denote the extension of the nontrivial automorphism of $L$ over $k$ to $E$ over $L$ by $\sigma$, which is the identity automorphism when restricted to $F$. 
Sometimes we write $\bar x=\sigma(x)$ for $x\in E$. 
Let us consider $E$ as a vector space over $L$. 
We consider the trace hermitian form on $E$ defined as follows: 
\begin{eqnarray*}
tr\colon  E\times E &\longrightarrow& L\\
tr(x,y)&=& tr_{E/L}(x\bar y).
\end{eqnarray*}
The restriction of this form to $F$ is the trace form $ tr\colon F\times F \longrightarrow k$, given by $tr(x,y)=tr_{F/k}(xy).$ 
We choose an orthogonal basis of $F$ over $k$, say $\{f_1, f_2, f_3\}$, with respect to the trace form, and extend it to a basis of $E/L$. 
Then the bilinear form $tr$ with respect to this basis has diagonalization $<1,2,2>$ (Section 18.31, \cite{kmrt}). 
We have $disc(tr)=4\in N_{L/k}(L^*)$. 
Hence $(E,tr)$ is a rank $3$ hermitian space over $L$ with trivial discriminant and $SU(E,tr)$ is isomorphic to $SU(H)$ where 
$H=diag(1,2,2)$.  
We choose an element $1\neq a \in T^1-L^1$, where $T^1=\{x\in E\mid x\bar x=1, N_{E/L}(x)=1\}$.  
Let us consider the left homothety map, 
\begin{eqnarray*}
l_a\colon  E &\longrightarrow& E\\
l_a(x)&=&ax
\end{eqnarray*}
Since $a\in T^1-L^1$, the characteristic polynomial $\chi(X)$ of $l_a$ is the minimal polynomial of $a$ over $L$, which is irreducible of degree $3$ over $L$.
 Next we prove that $l_a\in SU(E,tr)$. 
This is so since, 
$$
tr(l_a(x),l_a(y))=tr(ax,ay)=tr_{E/L}(ax\bar a\bar y)=tr_{E/L}(x\bar y)=tr(x,y).
$$
Let $S=(s_{ij})$ denote the matrix of $l_a$ with respect to the chosen basis $\{f_1, f_2, f_3\}$ of $F$ over $k$. 
Then the matrix of $l_{\bar a}$ is $\bar S=(\bar {s_{ij}})$. 
Also, since $a\bar a=1$, we have $S\bar S=1$. 
Thus we have a matrix $S$ in $SU(H)$, for $H=<1,2,2>$, satisfying the conditions of Proposition~\ref{counsun}. 
 
Now, let $L=k(\gamma)$ with $\gamma^2=c\in k^*$. 
We write $Q=k\oplus F$. 
Since $(F,tr)$ is a quadratic space with trivial discriminant, we can define a quaternionic multiplication on $Q$ (Prop. \ref{qb}), denote its norm by $N_Q$. 
We double $Q$ with $\gamma^2=c\in k^*$ to get an octonion algebra $\C=Q\oplus Q$ with multiplication,
$$
(x,y)(u,v)=(xu+c\bar vy,vx+y\bar u)
$$ 
and the norm $N((x,y))=N_Q(x)-cN_Q(y)$. 
We choose a basis $\{1,a,b,ab\}$ of $Q$, orthogonal for $N_Q$, so that $N_Q$ has diagonalization $<1,1,2,2>$ with respect to this basis. 
This gives a basis $\{(1,0),(a,0),(b,0),(ab,0),(0,1),(0,a),(0,b),(0,ab)\}$ of $\C$ and the diagonalization of $N$ with respect to this basis is $<1,1,2,2,-c,-c,-2c,-2c>$. 
We observe that the subalgebra $k\oplus k\subset \C$ is isomorphic to $L$ and $L^{\perp}=F\times F$ is a $3$ dimensional vector space over $L$ with hermitian form $<1,2,2>$. 
Hence $SU(L^{\perp},h)$, with respect to the basis $\{(a,0),(b,0),(ab,0)\}$ of $L^{\perp}$, is $SU(H)$ for $H=<1,2,2>$. 
Hence, from the discussion in previous paragraph, we have an element of required type in $SU(L^{\perp},h)$.  

By Corollary~\ref{lne}, $L^1/N(E^1)$ is nontrivial. 
It follows from Proposition~\ref{neum} and Proposition~\ref{counsun} that not all (semisimple) elements in $\Aut(\C)$, 
which are contained in the subgroup $SU(E,tr)$, are real. \ \ \ \ $\qed$
\begin{corollary}
Let $k$ be a totally real number field. 
Then there exists an octonion division algebra $\C$ over $k$ such that not every element in $\Aut(\C)$ is real. 
Hence there exist (semisimple) elements in $\Aut(\C)$, which are the product of three involutions but not the product of two involutions.
\end{corollary}
\noindent{\bf Proof : } We recall from Lemma~\ref{ohd} that if the $k$-quadratic form $q_B$, corresponding to the bilinear form 
$B\colon E\times E  \longrightarrow  k$, defined by $B(x,y)=tr_{E/L}(x\bar y)+tr_{E/L}(\bar x y)$, is anisotropic then the 
octonion algebra $\C$, as constructed in the proof of the above proposition, is a division algebra. 
In case when $k$ is a totally real number field and $L=k(i)$, the diagonalization of $q_B$ is $<1,2,2,1,2,2>$, which is clearly 
anisotropic over $k$. \ \ \ \ $\qed$

\noindent{\bf Remark: } We note that the quadratic form $q_B$ as above, can be isotropic for imaginary quadratic number fields.  
For example if $k=\Q(\sqrt {-2})$, $q_B$ has diagonalization $<1,-1,-1,-c,c,c>$, which is isotropic. 
Hence the octonion algebra $\C$ in this case is split. 
Therefore, indecomposable tori in subgroups $SU(V,h)\subset \Aut(\C)$ exist in all situations, whether $\C$ is division or not. 
And in either case, there are nonreal elements.

\subsection{$SL(3)\subset G$}\label{split}
Let us assume now that $L\cong k \times k$. 
We have seen in section 3 that $G(\C/L)\cong SL(U)\cong SL(3)$. 
Let $t_0$ be an element in $G(\C/L)$ and denote its image in $SL(3)$ by $A$. 
 We assume that the characteristic polynomial of $A\in SL(3)$ is irreducible over $k$. 
In this case, the characteristic polynomial equals the minimal polynomial of $A$. 

\begin{lemma}
Let the notations be fixed as above. 
Let $t_0$ be an element in $G(\C/L)$ and its image in $SL(3)$ be denoted by $A$. 
Let the characteristic polynomial of $A$ be irreducible over $k$. 
Suppose that $\exists h\in G=\Aut(\C)$, such that $ht_0h^{-1} =t_0^{-1} $. 
Then $h(L) = L$. 
\end{lemma}
\noindent
{\bf Proof :} 
Suppose that $h(L)\not\subset L$. 
Then we claim that $\exists x \in L\cap \C_0$ such that $h(x)\not\in L$. 
For this, let $y\in L$ be such that $h(y)\not\in L$. 
Let $x=y-\frac{1}{2}tr(y)1$. 
Then $tr(x)=0$ and if $h(x)\in L$ then $h(y)\in L$, a contradiction. 
Hence we have $x\in L\cap \C_0$ with $h(x)\not\in L$. 
Also since $t_0(x)=x$, we have, 
$$
t_0(h(x))=ht_0^{-1}(x)=h(x).
$$ 
Therefore, $t_0$ fixes $h(x)\in\C_0$ and hence fixes a two dimensional subspace $span\{x,h(x)\}$ pointwise, which is contained in $\C_0\subset \C$. 
Hence the characteristic polynomial of $t_0$ on $\C_0$ has a degree $2$ factor. 
But the characteristic polynomial $f(X)$ of $t_0$ on $\C_0$ has the factorization 
$$ 
f(X)=(X-1)\chi(X)\chi^*(X),
$$ 
where $\chi(X)$ is the characteristic polynomial of $t_0$ on the $3$ dimensional $k$-subspace $U$ of $\C_0$ and $\chi^*(X)$ is its dual polynomial (see Sec. 3,~\cite{w1}). 
Since $\chi(X)$ is irreducible by hypothesis, this leads to a contradiction. 
Hence any $h\in \Aut(\C)$, conjugating $t_0$ to $t_0^{-1}$ in $G$, leaves $L$ invariant.
\  \  \  \ $\qed$

\begin{theorem}\label{splitconjugate}
With notations fixed as above, let $A$ be the matrix of $t_0$ in $SL(3)$ with irreducible characteristic polynomial. 
Then $t_0$ is conjugate to $t_0^{-1}$ in $G=\Aut(\C)$, if and only if $A$ is conjugate to $\tr A$ in $SL(3)$. 
\end{theorem}
\noindent {\bf Proof : } Let $h\in G$ be such that $ht_0h^{-1}=t_0^{-1}$. 
In view of the lemma above, we have $h(L) = L$. 
We may, without loss of generality (up to conjugacy by an automorphism), assume that 
$$
\C =  \left \{ \left (\begin{array}{cc} \alpha &v \\  w&\beta \\ \end{array}  \right) \mid 
\alpha ,\beta \in k  ; v,w \in k^3 \right\}  \ \rm{with}\  L = \left \{ \left (\begin{array}{cc} \alpha & 0 
\\  0 &\beta \\ \end{array}  \right) \mid \alpha ,\beta \in k   \right\} 
$$ 
By Proposition~\ref{gl2}, $h$ belongs to $G(\C/L)\rtimes H$. 
Clearly $h$ does not belong to $G(\C/L) $, for if so, we can conjugate $t_0$ to $t_0^{-1}$ in $SL(U)$, which implies in particular that the characteristic polynomial $\chi(X)$ of $t_0$ on $U$ is reducible, a contradiction. 
Hence $h = g\rho$ for some $g \in G(\C/L)$. 
Let $A$ denote the matrix of $t_0$ on $U$ in $SL(3)$ and $B$ that of $g$. 
Then, a direct computation gives,
$$ 
h t_0 h^{-1} \left (\begin{array}{cc} \alpha &v \\  w&\beta  \end{array}  \right)
 =  \left (\begin{array}{cc} \alpha & B{\tr A}^{-1}B^{-1}v  \\ {\tra B}^{-1}A\tra B w &\beta \end{array}  \right),
$$
and
$$
t_0^{-1} \left (\begin{array}{cc} \alpha &v \\  w&\beta \end{array}  \right) 
=  \left (\begin{array}{cc} \alpha & A^{-1}v \\  \tr A w&\beta \end{array}  \right).
$$ 
Therefore,
$$
ht_0 h^{-1} = t_0^{-1} \Leftrightarrow A =B \tr A B^{-1}.
$$ 
Hence, $t_0$ is conjugate to $t_0^{-1}$ in $G(k)$ if and only if $A$ is conjugate to $\tra A$ in $SL(3)$.
\  \  \  \ $\qed$

We now derive a necessary and sufficient condition that a matrix $A$ in $SL(3)$, with irreducible characteristic polynomial, be conjugate 
to $\tr A$ in $SL(3)$. 
We have, more generally,
\begin{theorem}\label{trcon}
Let $A$ be a matrix in $SL(n)$ with characteristic polynomial $\chi_A(X)$ irreducible. 
Let $E=k[X]/\chi_A(X)\cong k[A]$ be the field extension of $k$ of degree $n$ given by $\chi_A(X)$.  
Then $A$ is conjugate to $\tr A $ in $SL(n)$, if and only if, for every $T\in GL(n)$ with $TAT^{-1} = \tr A$, $\det (T)$ is a norm from $E$.
\end{theorem}
\noindent {\bf Proof : } Fix a $T_0\in GL(n)$ such that $T_0AT_0^{-1} = \tr A$ and define a map, 
\begin{eqnarray*}
\{T\in M_n(k)\mid TA = \tr AT \} &\longrightarrow & k[A] \\
T&\mapsto& T_0^{-1}T
\end{eqnarray*}
This map is an isomorphism of vector spaces. 
Since if $T\in M_n(k)$ is such that $TA=\tr AT$ then $T_0^{-1}T$ belongs to $\mathcal Z(A)~(= k[A]$, as the characteristic polynomial of 
$A$ is the same as its minimal polynomial). 
To prove the assertion, suppose $T_0\in SL(n)$ conjugates $A$ to $\tr A$. 
But with the above bijection, $T_0^{-1}T = p(A)$  for some $p(A) \in k[A],~p(X)\in k[X]$. 
Hence $\det(p(A)) = \det(T)$, i.e. $\det T$ is a norm from $E$.

Conversely suppose there exists $T\in GL(n)$ with $TAT^{-1} = \tr A$ and $\det (T)$ is a norm from $E$. 
Then there exists $p(X) \in k[X]$ such that $\det(p(A)) = \det (T)^{-1}$. Thus $\det(T p(A)) =1$ and \\$(Tp(A))A(p(A)^{-1}T^{-1})=TAT^{-1}=\tr A.$ \ \ \ \ $\qed$

In the case under discussion, $A\in SL(3)$ has irreducible characteristic polynomial. 
Hence, $E \cong k[A]\cong \mathcal Z_{M_3(k)}(A)$  is a cubic field extension of $k$ . 
We combine the previous two theorems to get,
\begin{corollary}\label{example}
Let $A$ be a matrix in $SL(3)$ with irreducible characteristic polynomial. 
With notations as above, suppose $k^*/N(E^*)$ is trivial. Then $A$ can be conjugated to $\tr A$ in $SL(3)$ and hence $t_0$ can be conjugated to $t_0^{-1}$ in $\Aut(\C)$. 
\end{corollary}

If $k$ a $C_1$ field (e.g., a finite field) or $k$ does not admit cubic field extensions (e.g., $k$ real closed, algebraically closed), 
the above criterion is satisfied automatically. Hence every element in $G(\C/L)$, for $L=k\times k$, with irreducible characteristic polynomial over $k$, is conjugate to its inverse in $G(k)$. In particular, combining this with Theorem \ref{red}, we see that every semisimple element in $G(k)$ is real.  

We shall give a cohomological proof of reality for $G_2$ over fields $k$ with $cd(k)\leq 1$ (see the remarks later in this section).
\begin{lemma}\label{consln}
Let $A$ be a matrix in $SL(n)$ with irreducible characteristic polynomial. 
Then $A$ is conjugate to $\tr A$ in $SL(n)$ if and only if $A$ is a product of two symmetric matrices in $SL(n)$.
\end{lemma}
\noindent {\bf Proof : } Any matrix conjugating $A$ to $\tr A$ is necessarily symmetric (Th. 2,~\cite{tz}). 
Let $S$ be a symmetric matrix which conjugates $A$ to $\tr A$ in $SL(n)$, i.e., $SAS^{-1}=\tr A$. 
Let $B=SA=\tr A S$. 
Then $B$ is symmetric and belongs to $SL(n)$. Hence $A=S^{-1}B$ is a product of two symmetric 
matrices in $SL(n)$. 
Conversely, let $A$ be a product of two symmetric matrices from $SL(n)$, say $A=S_1S_2$. 
Then $S_2$ conjugates $A $ to $\tr A$. \ \ \ \ $\qed$

We need the following result from (\cite{w1}), (cf. also~\cite{l}),
\begin{proposition}\label{insym}
Let $\C$ be a (split) Cayley algebra over a field $k$ of characteristic not $2$. Let $L$ be a split two-dimensional subalgebra of $\C$. 
An element $\eta\in G(\C/L)$ is a product of two involutory automorphisms if and only if the corresponding matrix in $SL(3)$ can be 
decomposed into a product of two symmetric matrices in $SL(3)$.
\end{proposition}
We have,
\begin{theorem}\label{counterexample}
Let $t_0$ be an element in $G(\C/L)$, with notations as in this section. 
Let us assume that the  matrix $A$ of $t_0$ in $SL(3)$ has irreducible characteristic polynomial. 
Then, $t_0$ can be conjugated to $t_0^{-1}$ in $G=Aut(\C)$, if and only if $t_0$ is a product of two involutions in $G(k)$.
\end{theorem}
\noindent {\bf Proof : }  The element $t_0$ can be conjugated to $t_0^{-1}$ in $G$ if and only if, $A$ can be conjugated to $\tr A$ in $SL(3)$ (Theorem~\ref{splitconjugate}). 
This is if and only if, $A$ is a product of two symmetric matrices in $SL(3)$ (Lemma~\ref{consln}). 
By Proposition~\ref{insym}, this is if and only if $t_0$ is a product of two involutions in $\Aut(\C)$. \ \ \ \ $\qed$

In view of these results, to produce an example of a semisimple element of $G=\Aut(\C)$ that is not conjugate to its inverse in $\Aut(\C)$, 
we need to produce a semisimple element which is a product of three involutions but not a product of two involutions.
 We shall show that, for the split form $G$ of $G_2$ over $k=\Q$ or $k=\Q_p$, there are semisimple elements in $G(k)$ which are not conjugate to their inverses in $G(k)$. 
We shall end this section by exhibiting explicit elements in $G_2$ over a finite field, which are not real. 
These necessarily are not semisimple or unipotent (see the remarks at the end of this section).  
We adapt a slight variant of an example in (\cite{w1},\cite{l}) for our purpose, there the issue is bireflectionality of $G_2$. 
\begin{lemma}\label{cexample}
Let $k$ be a field and let $S$ be a symmetric matrix in $SL(3)$ whose characteristic polynomial $p(X)$ is irreducible over $k$. 
Let $E = k[X]/<p(X)>$, the degree three field extension of $k$ given by $p(X)$. 
Further suppose that $k^*/N(E^*)$ is not trivial. 
Then there exists a matrix in $SL(3)$, with characteristic polynomial $p(X)$, which is not a product of two symmetric matrices in $SL(3)$.  
\end{lemma}
\noindent{\bf Proof :}
Let $b \in k^*$ such that $b^2\not\in N(E^*)$. 
Consider $D = diag(b,1,1)$, a diagonal matrix and put $A=DSD^{-1}$. 
Then $A\in SL(3)$. 
We claim that $A$ is not a product of two symmetric matrices from $SL(3)$. 
Assume the contrary. 
Suppose $A=DSD^{-1}=S_1S_2$ where $S_1,S_2\in SL(3)$ and symmetric. Then 
$$
A=DSD^{-1}=(DSD)D^{-2}=S_1S_2.
$$ 
Let $T_1=DSD,~T_2=D^{-2}$. 
Then $ ^{t}T_i=T_i,~i=1,2$ and $A=T_1T_2=S_1S_2$. 
Therefore, 
$$
\tr A=T_2T_1=T_2AT_2^{-1}=S_2S_1=S_2AS_2^{-1}.
$$ 
Since the characteristic polynomial of $A$ is irreducible, by the proof of Theorem~\ref{trcon}, 
$D^2S_2=T_2^{-1}S_2\in\mathcal Z(A)=k[A]\cong E$.  
Which implies $S_2=D^{-2}f(A)$ for some polynomial $f(X)\in k[X]$. 
Taking determinants, we get 
$$
1=\det S_2=\det D^{-2}\det (f(A)),
$$ 
i.e., $b^{2}=\det(f(A))\in N(E^*) $, contradicting the choice of $b\in k$. 
Hence $A$ can not be written as a product of two symmetric matrices from $SL(3)$. \ \ \ \ $\qed$

\noindent{\bf Remark: } In view of Theorem~\ref{counterexample} and its proof, the element $A$ corresponds to an element in $\Aut (\C)$ which can not be conjugated to its inverse. 
If we choose the matrix $S$ as in the statement of the lemma above, to have separable characteristic polynomial, the matrix $A$, as constructed in the proof, corresponds to a semisimple element in an indecomposable torus contained in $SL(3)\subset G=\Aut(\C)$, which is not real.

\begin{theorem}\label{exsln}
Let $G$ be a split group of type $G_2$ over $k=\Q$ or $\Q_p$. 
Then there exists a semisimple element in $G_2(k)$ which is not conjugate to its inverse.
\end{theorem}
\noindent
{\bf Proof : }
{\bf Reality over $\Q_p$ : }
Let $k=\Q_p,~p\neq 2$. 
Let $p(X)$ be an irreducible monic polynomial of degree $n$, with coefficients in $\Q_p$. 
By a theorem of Bender (\cite{be1}), there exists a symmetric matrix with $p(X)$ as its characteristic polynomial, 
if and only if, for the field extension $E=\Q_p[X]/(p(X))$, there exists $\alpha$ in $E^*$ 
such that $(-1)^{\frac{n(n-1)}{2}}N(\alpha)$ belongs to $(\Q_p^*)^2$. 
We choose $E$ as the (unique) unramified extension of $\Q_p$ of degree $3$. Then, $E$ is a cyclic extension of $\Q_p$. 
We choose $\beta\in E^*,~N(\beta)=1$ so that  $E=\Q_p(\beta)$. 
Let $p(X)$ be the minimal polynomial of $\beta$ over $\Q_p$. 
Then, applying Bender's result, there is a symmetric matrix $A$ over $\Q_p$, with characteristic polynomial $p(X)$. 
Since $N(\beta)=1,$ $A$ belongs to $SL(3,\Q_p)$. 
We have, $\Q_p^*/N(E^*)\cong \Z/3\Z$ (see Sec. 17.9, \cite{p}),  
hence $(\Q_p^*)^2\not\subset N(E)$. 
Therefore we are done by Lemma~\ref{cexample}, combined with Proposition~\ref{insym} and Theorem~\ref{counterexample}.\ \ \ \ $\qed$

This example shows that there exist semisimple elements in $G=\Aut(\C)$ over $k=\Q_p$, which are not a product of two involutions and hence must be product of three involutions, by (\cite{w1}). 
In particular, reality for $G_2$ fails over $\Q_p$ (Theorem~\ref{maintheorem}).

\noindent
{\bf Reality over $\Q$ :}
A polynomial $p(X)\in K[X]$ is called $K$-real if  every real closure of $K$ contains the splitting field of $p(X)$ over $K$. 
Bender (Th. 1, \cite{be2}) proves that whenever we have $K$, an algebraic number field, and $p(X)$ a monic $K$-polynomial with an odd degree factor over $K$, then $p(X)$ is $K$-real if and only if it is the characteristic polynomial of a symmetric $K$-matrix. 

Let $p(X)=X^3-3X-1$. 
Then all roots of this polynomial are real but not rational. 
This polynomial is therefore irreducible over $\Q$ and by Bender's theorem stated above, $p(X)$ is the characteristic polynomial of a symmetric matrix. 
Note that $K=\Q[X]/<p(X)>$ is a degree $3$ cyclic extension of $\Q$.

We recall that for a cyclic field extension $K$ of $k$, the relative Brauer group $ B(K/k)\cong k^*/N_{K/k}(K^*)$ (ref. Sec. 15.1, Prop. b, \cite{p}). 
It is known that  if $K/k$ is a nontrivial extension of global fields, then $B(K/k)$ is infinite (ref. Cor. 4, \cite{fks}). 
Therefore, for $K$ chosen as above, $\Q^*/N(K^*)$ is not trivial. 
Hence all conditions required by Lemma~\ref{cexample} are satisfied by the polynomial $p(X)$ and we get a semisimple element $t_0 \in G_2(\Q)$ which is not conjugate to its inverse, using Lemma~\ref{cexample}, Proposition~\ref{insym} and Theorem~\ref{counterexample}. \ \ \ \ \ $\qed$

\noindent
{\bf Reality over $\F_q$ :} Let $k=\F_q$ be a finite field. 
We have shown (Th.~\ref{maintheorem}) that semisimple elements and unipotent elements in $G(k)$ are real in $G(k)$. 
We now construct an element in $G(k)$ which is not conjugate to its inverse. 
 Let $\C$ be the split octonion algebra over $k$, assume that the characteristic of $k$ is not $2$ or $3$. 
We use the matrix model for the split octonions, as introduced in the section  \ref{octonions}. 
Let $L$ be the split diagonal subalgebra of $\C$. 
We assume that $k$ contains primitive third roots of unity. 
We have, $G(\C/L)\cong SL(3)$. 
Let $\omega$ be a primitive third root of unity in $k$. 
Let 
$$
A=\left(\begin{array}{ccc} \omega & -1 & 0\\ 0 & \omega & 1 \\ 0 & 0 & \omega\end{array}\right).
$$ 
Then $A\in SL(3)$ and the minimal polynomial ($=$characteristic polynomial) of $A$ is $p(X)=(X-\omega)^3$.  
Let $b\in k$ be such that the polynomial $X^3-b^2$ is irreducible over $k$ (this is possible due to characteristic assumptions). 
Let $D=diag(b,1,1)$ and $B=DAD^{-1}$. 
Then $B\in SL(3)$ and has the same minimal polynomial as $A$. 
Note that $B$ is neither semisimple, nor unipotent. 
Let $t\in G(\C/L)$ be the automorphism of $\C$ corresponding to $B$. 
It is clear that the fixed point subalgebra of $t$ is precisely $L$.
\begin{theorem}\label{finitecounter}
The element $t\in G(\C/L)$ as above, is not real.
\end{theorem}
\noindent{\bf Proof:} If not, suppose for $h\in G(k)$, $hth^{-1}=t^{-1}$. 
Then, since $t$ fixes precisely $L$ pointwise, we have $h(L)=L$. 
Therefore $h\in G(\C,L)\cong G(\C/L)\rtimes H$, where $H=<\rho>$ is as in Prop.~\ref{gl2}. 
If $h\in G(\C/L)$, conjugacy of $t$ and $t^{-1}$ by $h$ would imply conjugacy of $B$ and $B^{-1}$ in $SL(3)$. 
But this can not be, since $\omega$ is the only root of $p(X)$. 
Thus $h=g\rho$ for $g\in G(\C/L)$. 
Now, by exactly the same calculation as in the proof of Theorem~\ref{splitconjugate}, conjugacy of $t$ and $t^{-1}$ in $G(k)$ is equivalent to 
conjugacy of $B$ and $\tra B$ in $SL(3)$. 
Let $CBC^{-1}=\tra B$ with $C\in SL(3)$. 
Let 
$$
T=\left(\begin{array}{ccc}0 & 0 & 1\\ 0 & -1 & 0\\ 1 & 0 & 0\end{array}\right).
$$ 
Then $T\in SL(3)$ is symmetric and $TAT^{-1}=\tr A$. 
Hence $A$ is a product of two symmetric matrices in $SL(3)$, say $A=T_1T_2$ with $T_i\in SL(3)$, symmetric (see the proof of Lemma~\ref{consln}). 
But $CBC^{-1}=\tra B$ gives $(DCD)A=\tr A (DCD)$. 
Therefore, by an argument used in the proof of Theorem~\ref{trcon}, using the fact that the characteristic polynomial is equal to the minimal polynomial of $A$, 
we have, $DCD=T_2f(A)$ for some polynomial $f\in k[X]$. 
Taking determinants, we get $b^2=\det (f(A))=f(\omega)^3$. 
But this contradicts the choice of $b$. Hence $t$ is not real. \ \ \ \ $\qed$

A similar construction can be done for the subgroup $SU(V,h)\subset G$. 
We continue to assume that $k$ is a finite field with characteristic different from $2, 3$. 
We first note that the (split) octonion algebra contains all quadratic extensions of $k$. 
We assume that $2$ is a square in $k$ and that $k$ contains no primitive cube roots of unity. 
Let $L$ be a quadratic extension of $k$ containing a primitive cube root of unity $\omega$. 
Let $b\in L$ with $N_{L/k}(b)=1$ such that the polynomial $X^3-b^2$ is irreducible over $L$. 
Let $\alpha\in L$ with $N_{L/k}(\alpha)=-1$. 
Let 
$$
A=\left(\begin{array}{ccc}\omega+\frac{1}{4} & \frac{1}{2} & -\frac{1}{4}\alpha\\ -\frac{1}{2}\omega^2 & \omega & \frac{1}{2}\alpha\omega^2\\ -\frac{1}{4}\overline{\alpha} & -\frac{1}{2}\overline{\alpha} & \omega-\frac{1}{4}\end{array}\right),
$$
then $A\in SU(3)$ and the minimal polynomial ($=$ characteristic polynomial) of $A$ over $L$ is $(X-\omega)^3$. 
Let $F$ be a cubic extension of $k$ and $E=F.L$.
 Then $E$ is a cyclic extension of $L$ and we have the trace hermitian form as defined in Prop.~\ref{exsuh}, on $E$. 
We fix an orthogonal basis for $F$ over $k$ for the trace bilinear form and extend it to a basis of $E$ over $L$. 
Then the trace hermitian form has diagonalization $<1,1,1>$. 
We construct $\C=L\oplus E$ with respect to the hermitian space $(E,tr)$, as in Section $3$. 
Then $SU(L^{\perp},h)\cong SU(3)$. 
Let $D=diag(b,1,1)$ and $B=DAD^{-1}$. 
Then $B\in SU(3)$ and has the same minimal polynomial as $A$. 
Note that $B$ is neither semisimple, nor unipotent. 
Let $t$ denote the automorphism of $\C$ corresponding to $B$. 
Then the fixed point subalgebra of $t$ in $\C$ is precisely $L$. 
We have, 
\begin{theorem}\label{finitecounter'}
The element $t\in G(\C/L)$ as above, is not real.
\end{theorem}
\noindent{\bf Proof:} Suppose $t$ is real in $G(k)$. 
Then there is $h\in G(k)$ such that $hth^{-1}=t^{-1}$. 
Since the fixed point subalgebra of $t$ is $L$, we have $h(L)=L$. 
Thus, by Proposition~\ref{gl21}, $h\in G(\C,L)\cong G(\C/L)\rtimes H$, where $H=<\rho>$ is as in Proposition~\ref{gl21}. 
If $h\in G(\C/L)$, then $B$ and $B^{-1}$ would be conjugate in $SU(3)$, but that can not be since $\omega$ is the only eigenvalue for $B$. 
Hence $h=g\rho$ for $g\in G(\C/L)$. 
Then, conjugacy of $t$ and $t^{-1}$ in $G(k)$ is equivalent to conjugacy of $\overline{B}$ and $B^{-1}$ in $SU(3)$, by the same calculation as in the proof of Th. 6.5. 
By Lemma~\ref{suco1}, this is if and only if $B=B_1B_2$ with $B_i\in SU(3)$ and $\overline{B_i}B_i=1$. 
But then $B=DAD^{-1}=B_1B_2$ and hence $A=(D^{-1}B_1D^{-1})(DB_2D)=A_1A_2$, say. 
Then $A_i\in U(3)$ and $\overline{A_i}A_i=1$. 
Let $C\in SU(3)$ be such that $C\overline{B}C^{-1}=B^{-1}$. 
Then $C\overline{DAD^{-1}}C^{-1}=DA^{-1}D^{-1}$. 
This gives, $(D^{-1}CD^{-1})\overline{A}(DC^{-1}D)=A^{-1}$. 
Hence, by Lemma~\ref{suco1}, $A=T_1T_2$ with $T_i\in SU(3),~\overline{T_i}T_i=1$. 
Therefore, by a similar argument as in the remark following Corollary~\ref{suobs}, we must have, $T_1A_1^{-1}=f(A)$ for a polynomial $f(X)\in L[X]$. 
Taking determinants, we get $b^{-2}=f(\omega)^3$, contradicting the choice of $b$. 
Therefore $t$ is not real in $G(k)$.  \ \ \ \ $\qed$

\paragraph*{\bf Remarks :}
\begin{enumerate}
\item[1.] Our results in fact show that if a semisimple element in $G(k)$, for a group $G$ of type $G_2$, is conjugate to its inverse in $G(k)$, the conjugating element can be chosen to be an involution. 
The same is true for unipotents (these are always conjugate to their inverses).
\item[2.] One can give a simple cohomological proof of reality for $G_2$ over $k$ with $cd(k)\leq 1$. 
Recall that $cd(k)\leq 1$, if and only if for every algebraic extension $K$ of $k,~Br(K)=0$ (Chap. 3, Prop. 5, \cite{se}).
Let $g\in G(k)$ be semisimple and  $T$ be a maximal $k$-torus of $G$ containing $g$ (cf. Corollary 13.3.8, \cite{sp}). 
Let $N(T)$ be the normalizer of $T$ in $G$ and $W=N(T)/T$ the Weyl group of $G$ relative to $T$. 
We have the exact sequence of groups,
$$
1\rightarrow T\rightarrow N(T)\rightarrow W\rightarrow 1.
$$ 
The corresponding Galois cohomology sequence is,   
$$
1\rightarrow T(k)\rightarrow N(T)(k)\rightarrow W(k)\rightarrow H^1(k,T)\rightarrow\cdots
$$
Now, if $cd(k)\leq 1$, by Steinberg's theorem (See~\cite{s}), $H^1(k,T)=0$. 
Hence the last map above is surjective homomorphism of groups. 
Therefore the longest element $w_0$ in the Weyl group $W$ of $G_2$, which acts as $-1$ on the set of positive roots with respect to $T$ (\cite{b}, Plate IX), lifts to an element 
$h$ of $N(T)(k)$. 
Hence, over $k_s$, we have $hth^{-1}=t^{-1}$ for all $t\in T$. 
But $h\in G(k)$, hence the conjugacy holds over $k$ 
itself. 
Using Theorem 6.3, we get the following interesting result,
\begin{theorem}\label{cdk}
Let $cd(k)\leq 1$ and $G$ be a group of type $G_2$ over $k$. 
Then every semisimple element in $G(k)$ is a product of two involutions in 
$G(k)$.
\end{theorem}   
\noindent
\item[3.] {\bf The obstruction :} From our results, we see that semisimple elements belonging to decomposable tori are always product of two involutions and hence real in $G(k)$. 
For semisimple elements belonging to an indecomposable maximal torus $T$, the obstruction to reality is measured by $L^1/N(\E^1)$, 
where $T\subset SU(V,h)\cong SU(\E,h^{(u)})$ is given by $T=\E^1$ and $\E$ is a cubic field extension of $L$. 
In the other case, when $T\subset SL(3)$, the obstruction is measured by $k^*/N(\F^*)$, where $\F$ is a cubic field extension of $k$. 
In both cases, the obstruction has a Brauer group interpretation. When $T\subset SL(3)\subset G$ is an indecomposable maximal torus, 
coming from a cyclic cubic field extension $\F$ of $k$, the obstruction to reality for elements in $T(k)$, is the relative Brauer group 
$B(\F/k)$. For an indecomposable torus $T\subset SU(\E,h^{u})\subset G$, where $\E$ is a cubic cyclic field extension of $L$, 
the obstruction is the quotient $B(\E/L)/\phi(B(\F/k)$, where $\F$ is the subfield of $\E$, fixed by the involution $\sigma$ on $\E$ and $\phi$ is the base change map $B(\F/k)\longrightarrow B(\E/L)$.
\end{enumerate}           
\noindent
{\bf Acknowledgments :}
We thank Dipendra Prasad for suggesting the problem and his generous help. 
We thank Huberta Lausch and E. Bender for making available their papers to us. 
The authors thank H. Petersson for very useful comments on the paper. 
We thank Surya Ramana and Shripad for many stimulating discussions. 
We are extremely grateful to the referee for suggestions, 
which improved the exposition tremendously. 
The second author thanks ICTP-Trieste for its hospitality during the summer of 2003.

\end{document}